\documentclass[12pt,reqno]{amsart}   

\usepackage{graphicx,subfigure,color}
\usepackage{latexsym,amsmath,amsfonts,amscd,amsthm,epstopdf,lineno}
\usepackage[T1]{fontenc}
\usepackage{caption,here}
\usepackage{textcomp}
\usepackage{tikz}
\usetikzlibrary{decorations.pathmorphing,decorations.pathreplacing}
\usepackage{siunitx}
\usepackage{algorithm}  
\usepackage{enumitem}   
\usepackage{fullpage}
\usepackage[colorlinks=true,citecolor=blue,urlcolor=blue]{hyperref}
\newtheorem{remark}{Remark}

\newcommand{\Cross}{$\mathbin{\tikz [x=1.4ex,y=1.4ex,line width=.2ex] \draw (0,0) -- (2,2) (0,2) -- (2,0);}$}%

\newcommand{\dt}{\Delta t}
\newcommand{\dx}{\Delta x}
\newcommand{\ainvs}{\frac{1}{\alpha^{2}}}

\begin{document}

\title{Method of lines transpose:  an efficient A-stable solver for wave propagation}
\author[M. Causley]{Matthew Causley}
\address{Mathematics Department, Kettering University, Flint, MI 48504}
\email{mcausley@kettering.edu}

\thanks{This work has been supported in part by AFOSR grants FA9550-11-1-0281, FA9550-12-1-0343 and FA9550-12-1-0455, NSF grant DMS-1115709, and MSU Foundation grant SPG-RG100059.}

\author[A. Christlieb]{Andrew Christlieb}
\address{Department of Mathematics, Michigan State University, East Lansing, MI 48824}
\email{andrewch@math.msu.edu}

\author[E. Wolf]{Eric Wolf}
\address{Department of Mathematics, Michigan State University, East Lansing, MI 48824}
\email{wolferi@math.msu.edu}


\date{}

\begin{abstract}
Building upon recent results obtained in \cite{Causley2013a,Causley2013,Causley_Christlieb_Cho}, we describe an efficient second order, A-stable scheme for solving the wave equation, based on the method of lines transpose (MOL$^T$), and the resulting semi-discrete (i.e. continuous in space) boundary value problem. In \cite{Causley2013a}, A-stable schemes of high order were derived, and in \cite{Causley_Christlieb_Cho} a high order, fast $\mathcal{O}(N)$ spatial solver was derived, which is matrix-free and is based on dimensional-splitting.

In this work, are interested in building a wave solver, and our main concern is the development of boundary conditions. We demonstrate all desired boundary conditions for a wave solver, including outflow boundary conditions, in 1D and 2D. The scheme works in a logically Cartesian fashion, and the boundary points are embedded into the regular mesh, without incurring stability restrictions, so that boundary conditions are imposed without any reduction in the order of accuracy.  We demonstrate how the embedded boundary approach works in the cases of Dirichlet and Neumann boundary conditions. Further, we develop outflow and periodic boundary conditions for the MOL$^T$ formulation.  Our solver is designed to couple with particle codes, and so special attention is also paid to the implementation of point sources, and soft sources which can be used to launch waves into waveguides.


\end{abstract}

\maketitle

\section{Introduction}

Numerical solutions of the wave equation have been an area of investigation for many decades. The wave equation is ubiquitous in the physical world, arising in acoustics, electromagnetics, and fluid dynamics. Our main interest is in electromagnetic wave propagation, where traditional finite difference methods such as the finite-difference time-domain (FDTD) method are often used to solve Maxwell's equations. When the classical Yee scheme is used, the Courant-Friedrichs-Lewy (CFL) stability criterion restricts the time step to scale with the smallest cells  in the domain. This becomes computationally prohibitive in a variety of interesting problems, such as electromagnetic scattering or waveguide design, where complex geometries need to be embedded in a Cartesian mesh, leading to small spatial cells near the boundaries. Alternatively, problems with multiple and disparate temporal or spatial scales, such as those presented by plasma simulations, require time steps and mesh spacings which are on the order of the shortest temporal and spatial scales. Due to the large number of charged particles in a typical simulation \cite{Birdsall1976}, and the high dimensionality of their distribution space, time steps are computationally prohibitive, and the amount of them must be minimized.

As a result, A-stable Maxwell solvers have been developed, which remove the CFL restriction and restore the ability of the user to define a time step which is based on the physical problem, rather than on its spatial discretization. Notable advances include the introduction of the FDTD-ADI \cite{Fornberg,Fornberga,Namiki2000,Smithe2009} algorithm, as well as several semi-implicit time split schemes. However, it remains difficult to preserve both accuracy, and A-stability when non-rectangular domains are required.

Alternatively, Maxwell's equations can be reformulated (e.g., using the scalar and vector potential formulation), so that each field component independently satisfies the second order wave equation. Method of lines (MOL) formulations of the wave equation are well-studied, and often lead to conditionally stable schemes. But when the discretization is first performed in time, following the method of lines transpose (MOL$^T$), the wave equation is found at discrete time levels by solving a semi-discrete boundary value problem. In \cite{Alpert2000, Alpert2002}, the MOL$^T$  is used to produce an exact integral formulation of the wave equation, where the solution is determined by a convolution over the domain of dependence against a space-time Green's function $G_d(x,t)$, in dimensions $d = 1,2, 3$. However, more commonly the semi-discrete solution is obtained, in which the modified Helmholtz equation must be solved, and the corresponding semi-discrete Green's function $G(x,\Delta t)$ is the Yukawa, or modified Helmholtz kernel. The resulting boundary integrals methods can then be solved at each time step using fast summation algorithms, such as tree-codes \cite{Barnes1986, Christlieb2004, Li2009, Lindsay2001}, or the fast multipole method (FMM) \cite{greengard1987fast, Cheng1999, Coifman1993, Li2006, Gimbutas2002, Li2009, Cheng2006}. These methods often scale as $O(N)$ or $O(N\log N)$, but require substantial storage or precomputing stages.

It is worth noting however that in one spatial dimension, the modified Helmholtz kernel is a simple decaying exponential. This fact has been combined with alternate dimension implicit (ADI) splitting to achieve an A-stable wave solver with computational complexity of $O(N \log N)$ \cite{Bruno2010, Lyon2010}, and $O(N)$ \cite{Causley2013, Causley2013a}, respectively. The scaling presented in \cite{Bruno2010, Lyon2010} is due to the Fourier continuation method, which makes use of the fast Fourier transform to compute a periodic extension of the the convolution integral, where the period is taken sufficiently longer than the domain, so that boundary conditions can be satisfied.

In our approach \cite{Causley2013a, Causley2013,Causley_Christlieb_Cho}, the one-dimensional convolution is instead computed strictly over the computational domain, using polynomial interpolation. In \cite{Causley2013}, the solution is proven to be A-stable, second order accurate, and that the matrix formed by the discrete convolution can be applied in $\mathcal{O}(N)$ operations. However, the discrete convolution matrix was formed using piecewise linear integration, and it was found that a Lax-type correction was required to ensure convergence of the scheme in the semi-discrete limit $\Delta t\to 0$, where $\Delta x$ is fixed. This issue was addressed in \cite{Causley_Christlieb_Cho}, where spatial discretization was extended to orders $M\geq 2$, on non-uniform grids while retaining  $\mathcal{O}(N)$ complexity. In \cite{Causley2013a}, the scheme was extended to higher orders in time using a novel approach, successive convolution. These higher order schemes were proven to be A-stable.

The purpose of this present work is to develop a robust embedded boundary approach for complex geometry for the MOL$^T$ formulation.  Because this paper is focused on developing a range of boundary conditions for the wave equation, we limit our attention to the 2nd order accurate solver. In addition to addressing standard closed (i.e., Dirichlet, Neumann and periodic) boundary conditions, we also develop an open, or outflow boundary condition which is suitable for our implicit solver.  The method is extended to higher spatial dimensions using the factorization developed in \cite{Causley2013a}, which is similar to but different from the traditional ADI splitting.   This means that solutions are computed line-by-line along dimensional sweeps leveraging the $O(N)$ 1D solver to construct a high dimensional $O(N$) implicit method.

The rest of this paper is laid out as follows. In section \ref{sec:molt}, we use the method of lines transpose (MOL$^T$) to reformulate the 2D wave equation as an implicit, semi-discrete modified Helmholtz equation. The Helmholtz operator is inverted analytically using dimensional splitting, and we recast the solution as a series of one-dimensional boundary integral equations. We specifically show how to obtain a fully discretized, second order accurate solution by performing spatial quadrature, using a fast $\mathcal{O}(N)$ convolution algorithm.

In section \ref{sec:artificial_dissipation}, we modify the time-centered scheme of \cite{Causley2013a} to introduce artificial dissipation. This will be used to stabilize the embedded boundary method in our 2D algorithm, presented in section \ref{sec:BC}. We demonstrate our solution to be fast, second order accurate, and A-stable, even on non-rectangular domains with numerical results presented in section \ref{sec:results}.  We conclude the body of the paper with several remarks in section \ref{sec:conclusion}.

We also include in Appendix \ref{sec:wave_summary}, a table summarizing the main algorithmic aspects of our wave solver. In Appendix \ref{sec:sources}, we show how point sources, which may represent charged plasma particles, or soft point sources launched into a waveguide, can be implemented.

\section{Integral solution using MOL$^T$}
\label{sec:molt}

We now develop a dimensionally split algorithm for solving the initial value problem
\begin{align}
	\nabla^2 u - \frac{1}{c^2}\frac{\partial^2 u}{\partial t^2} &= -S(\mathbf{x},t), \quad  \mathbf{x} \in \Omega, \quad t>0 \label{eqn:prob2} \\
	u(\mathbf{x},0) &= f(\mathbf{x}), \quad \mathbf{x} \in \Omega \nonumber \\
	u_t(\mathbf{x},0) &= g(\mathbf{x}), \quad \mathbf{x} \in \Omega \nonumber,
\end{align}
with consistent boundary conditions. We utilize the method of lines transpose (MOL$^T$) to perform a second order accurate temporal discretization, as was shown in \cite{Causley2013a}. In particular, the discretization is time-centered, and implicit in the spatial derivatives
\[
\left(1 - \frac{\nabla^2}{\alpha^2}\right) \left[u^{n+1}-2u^n+u^{n-1} +\beta^2 u^n\right] (x,y)=  \beta^2 \left(u^n + \frac{1}{\alpha^2} S^n \right),
\]
where
\[
\alpha = \frac{\beta}{c\Delta t}, \qquad \beta>0,
\]
and $\beta\leq 2$ is chosen to enforce A-stability \cite{Causley2013a}. Inversion of this operator leads to a boundary integral formulation, where the Green's function is the Yukawa potential, defined for 2D in terms of the modified Bessel function $K_0(r)$. However we will employ our previously developed \cite{Causley2013a,Causley2013} splitting algorithm, which has $\mathcal{O}(N)$ complexity
\begin{equation}
\label{eqn:Helmholtz_Equation_2D}
\mathcal{L}_x \mathcal{L}_y \left[u^{n+1}-2u^n+u^{n-1} +\beta^2 u^n\right] (x,y)=  \beta^2 \left(u^n + \frac{1}{\alpha^2} S^n \right).
\end{equation}
Here the subscripts denote the spatial component of univariate modified Helmholtz operators,
\begin{equation}
\label{eqn:HelmholtzL}
\mathcal{L}_x[u]: = \left(1 - \frac{\partial_{xx}}{\alpha^2}\right)u(x,y), \qquad	\mathcal{L}_y[u]: = \left(1 - \frac{\partial_{yy}}{\alpha^2}\right)u(x,y).
\end{equation}
The modified Helmholtz operators are formally inverted using the Green's function. We define convolution with this Green's function by the integral operator
\begin{equation}
\label{eqn:Iu}
I_x[u](x,y) := \frac{\alpha}{2}\int_a^b u(x',y)e^{-\alpha|x-x'|}dx', \quad a\leq x \leq b,
\end{equation}
so that
\begin{equation}
\label{eqn:L_Inverse}
\mathcal{L}_x^{-1}[u](x,y):= \underbrace{I_x[u](x)}_{\text{Particular Solution}}+ \underbrace{\vphantom{I_x[u](x)} A e^{-\alpha(x-a)} + B e^{-\alpha(b-x)}}_{\text{Homogeneous Solution}},
\end{equation}
where the coefficients $A$ and $B$ of the homogeneous solution are determined by applying boundary conditions. A similar definition holds for $\mathcal{L}_y^{-1}$.
Formally inverting both operators to the right hand side we find the explicit equation
\begin{equation}
\label{eqn:2D_explicit}
u^{n+1}-2u^n+u^{n-1} = -\beta^2 \mathcal{D}_{xy}[u^n] + \beta^2\mathcal{L}_{x}^{-1}\mathcal{L}_{y}^{-1}\left[\frac{1}{\alpha^2} S^n \right](x,y),
\end{equation}
where the multidimensional operator is now
\[
\mathcal{D}_{xy}[u]: = u - \mathcal{L}_x^{-1}\mathcal{L}_y^{-1}[u].
\]
Since the application of $\mathcal{L}_x^{-1}$ is done for fixed $y$, and vice versa for $\mathcal{L}_y^{-1}$, the operator $\mathcal{D}_{xy}$ can be constructed in a line-by-line fashion, similar to ADI algorithms. It was also proven in \cite{Causley2013a} that the scheme \eqref{eqn:2D_explicit} is A-stable, for $0< \beta \leq 2$.
\begin{remark}
	In the continuous case, $\mathcal{L}_x$ and $\mathcal{L}_y$ (and their inverses) commute, so $\mathcal{D}_{xy} = \mathcal{D}_{yx}$. However in practice the discretize operators will not commute, and some small spatial anisotropy is introduced. This can be controlled by applying both operators $\mathcal{D}_{xy}$ and $\mathcal{D}_{yx}$, and then averaging the result.
\end{remark}

\begin{remark}
	Similar to a traditional ADI formulation of the wave equation, this factorization produces a fourth order splitting error term,
	\[
	\mathcal{L}_x\mathcal{L}_y = \left(1-\frac{\partial_{xx}}{\alpha^2}\right)\left(1-\frac{\partial_{yy}}{\alpha^2}\right) = \left(1-\frac{\partial_{xx}+\partial_{yy}}{\alpha^2}+\frac{\partial_{xx}\partial_{yy}}{\alpha^4}\right),
	\]
	which can be compensated for by adding a term to the right hand side
	\begin{equation}
	\label{eqn:Helmholtz_Equation_2D_corrected}
	\mathcal{L}_x \mathcal{L}_y \left[u^{n+1}-2u^n+u^{n-1} +\beta^2 u^n\right] =  \beta^2 \left(u^n + \frac{1}{\alpha^2} S^n \right) +\beta^2\left(\mathcal{L}_x -1 \right)\left(\mathcal{L}_y -1 \right)[u^n].
	\end{equation}
	Note our use of the identity $(\mathcal{L}_x - 1)(\mathcal{L}_y - 1) = \partial_{xx}\partial_{yy}/\alpha^4$. Formally inverting both operators to the right hand side, the solution can be rearranged and found as
	\begin{equation}
	\label{eqn:2D_explicit_corrected}
	u^{n+1}-2u^n+u^{n-1} = - \beta^2 \mathcal{C}[u^n ](x,y) +\beta^2\mathcal{L}_{x}^{-1}\mathcal{L}_{y}^{-1}\left[ \frac{1}{\alpha^2} S^n \right](x,y),
	\end{equation}
	where the convolution operator $\mathcal{C}$ is the operator defined in \cite{Causley2013a} as
	\begin{equation}
	\label{eqn:C_def}
	\mathcal{C}:= \mathcal{L}_x^{-1} \mathcal{D}_y + \mathcal{L}_y^{-1}\mathcal{D}_x = \mathcal{D}_{xy} -\mathcal{D}_x\mathcal{D}_y.
	\end{equation}
	This form was used in \cite{Causley2013a} to achieve schemes of higher order through successive convolution, where removing splitting errors is of paramount concern.
\end{remark}

\subsection{Fast convolution algorithm for the integral solution}
\label{sec:Fast}
In \cite{Causley2013}, the particular solution \eqref{eqn:Iu} was discretized in space using the weighted midpoint and trapezoidal rules, which amounted to replacing $u$ with a piecewise constant and linear approximation, respectively. More recently \cite{Causley_Christlieb_Cho}, we have detailed the spatial discretization of $I_x$ to arbitrary order, while providing a means for its rapid evaluation. In this work we focus on the second order accurate implementation of this algorithm, and reiterate the relevant details.

This particular solution is first decomposed into a left and right oriented integral, split at $y=x$ so that
\begin{equation}
\label{eqn:ILR}
I[u](x) = I^L[u](x)+I^R[u](x),
\end{equation}
where
\[
I^L[u](x) = \frac{\alpha}{2} \int_a^x e^{-\alpha(x-y)}u(y)dy, \quad  I^R[u](x)= \frac{\alpha}{2}\int_{x}^b e^{-\alpha(y-x)}u(y)dy.
\]
These quantities can be updated locally, using exponential recursion
\begin{align}
	\label{eqn:IL_def}
	I^L[u](x) &= I^L[u](x-\delta_L) e^{-\alpha \delta_L}+ J^L[u](x), \quad J^L[u](x):= \frac{\alpha}{2} \int_{0}^{\delta_L} u(x-y) e^{-\alpha y}dy, \\
	\label{eqn:IR_def}
	I^R[u](x) &= I^R[u](x+\delta_R) e^{-\alpha \delta_R}+ J^R[u](x), \quad J^R[u](x):= \frac{\alpha}{2} \int_{0}^{\delta_R} u(x+y) e^{-\alpha y}dy.
\end{align}
The recursive updates \eqref{eqn:IL_def} and \eqref{eqn:IR_def} are exact (in space), and making $\delta_L$ and $\delta_R$ small (typically, $\delta_L =\delta_R = \Delta x$) effectively localizes the contribution of the integrals.

Based on these observations, we now outline the fast convolution algorithm. In \cite{Causley_Christlieb_Cho}, we derived a spatial quadrature of general order $M \geq 2$ for irregular grids. Here, we will utilize a second-order accurate method on a uniform grid, which we explicitly develop.

Consider the domain $(a,b)$ discretized by uniform grid points $x_{1} = a < x_{2} < \cdots < x_{N+1} = b$ of width $\Delta x = \frac{b-a}{N} = x_{j+1}-x_{j}$, $j=1,...,N$. Suppose, given a function $f$ compactly supported in $(a,b)$, we are to evaluate the convolution operator at each grid points, that is compute $I_{j} = I[f](x_{j}) = I^{L}[f](x_{j})+I^{R}[f](x_{j}) = I^{L}_{j}+I^{R}_{j}$, in $O(N)$ operations. We proceed by evaluating the local integrals $J^{L}_{j} = J^{L}[f](x_{j})$ and $J^{L}_{j} = J^{L}[f](x_{j})$, as in \eqref{eqn:IL_def} and \eqref{eqn:IR_def}. The local integrals may be evaluated with quadrature, or, if possible, analytically. A second-order accurate quadrature is given by 

\begin{align}
	J^{L}_{j} &\approx P f(x_{j})+Qf(x_{j-1})+R(f(x_{j+1})-2f(x_{j})+f(x_{j-1})) \\
	J^{R}_{j} &\approx P f(x_{j})+Qf(x_{j+1})+R(f(x_{j+1})-2f(x_{j})+f(x_{j-1}))
\end{align}
where defining $\nu = \alpha \dx$ and $d = e^{-\nu}$, the quadrature weights are given by
\begin{align}
	P&= 1-\frac{1-d}{\nu} \\
	Q&= -d+\frac{1-d}{\nu} \\
	R&= \frac{1-d}{\nu^{2}}-\frac{1+d}{2\nu}.
\end{align}

We summarize our fast method in Algorithm \ref{alg:fast_conv}.

\begin{algorithm}[htp]
	\caption{Fast convolution algorithm}
	\label{alg:fast_conv}
	\begin{enumerate}[itemsep=2mm]
		\item Compute $J^{L}_{j+1}$ and $J^{R}_{j}$ for $j=1,...,N$ via quadrature or analytical integration.
		
		\item Initialize $I^{L}_{1} = 0$ and $I^{R}_{N+1}=0$, and perform the exponential recursion,
		$I^{L}_{j+1} = J^{L}_{j+1} + e^{-\alpha \Delta x} I^{L}_{j}$ for $j=1,...,N$ and $I^{R}_{N-j+1} = J^{R}_{N-j+1} + e^{-\alpha \Delta x} I^{R}_{N-j+2}$ for $j=1,...,N$.

	\end{enumerate}
\end{algorithm}


\section{Artificial dissipation}
\label{sec:artificial_dissipation}

As will be discussed in section \ref{sec:embedded_neumann_bc}, it is necessary to include some artificial dissipation in the numerical scheme to maintain stability with embedded boundary methods for Neumann boundary conditions. We first present a version of the wave solver based on a backwards difference formula (BDF) time discretization, leading to what we call a diffusive scheme, which is dissipative. This method has a larger truncation error than a centered scheme, does not possess a means to tune the level of dissipation, and also has an implicit source term (at time level $n+1$), which is problematic for application in the context of particle-in-cell (PIC) methods for the simulation of plasmas. The second-order centered (dispersive) scheme given above is therefore preferable, but are non-dissipative in their original forms. We give a method for adding tunable artificial dissipation terms into the centered scheme, while maintaining A-stability.

\subsection{Diffusive wave solver}
We substitute the following backward difference formula (BDF) discretization:
\begin{align}
	u_{tt}^{n+1} &= \frac{2u^{n+1}-5u^n+4u^{n-1}-u^{n-2}}{\Delta t^2} - \frac{11\Delta t^2}{12} u_{tttt}(x,\eta)
\end{align}
into the wave equation $\frac{1}{c^{2}}u_{tt}-\nabla^{2} u = S(x,t)$.

Rearranging, defining $\alpha = \sqrt{2}/(c \Delta t)$ and dividing by $\alpha^{2}$ gives the semi-discrete scheme
\begin{align}
	\left(-\ainvs \Delta+1\right) u^{n+1} = \frac{1}{2}\left(5u^{n}-4u^{n-1}+u^{n-2}\right) + \ainvs S(x,t^{n+1}) + O(\dt^{4}).
\end{align}

This method is A-stable and dissipative, but does not possess a mechanism for tuning the dissipation, has an inconvenient implicit source term (at time level $n+1$), and typically has a larger truncation error compared to the centered scheme.

\subsection{Artificial dissipation in centered schemes}

\subsubsection{Artificial Dissipation in 1D}
We give a modified form of the centered version of the implicit wave solver with tunable artificial dissipation that retains the property of unconditional stability. We let $\epsilon$ denote a small artificial dissipation parameter, and $\mathcal{D}_x[u] = u-\mathcal{L}^{-1}[u] = u(x)-\frac{\alpha}{2}\int_{-\infty}^{\infty} e^{-\alpha|x-x'|}u(x') \, dx'$ be defined as usual.

Ignoring sources, we have the second order scheme with dissipation,

\begin{align}
	u^{n+1}-2u^{n}+u^{n-1} = -\beta^{2}\mathcal{D}_x[u^{n}]+ \epsilon \mathcal{D}_x^{2}[u^{n-1}].
\end{align}

We now prove the unconditional stability of this scheme for prescribed values of $\beta$. As in \cite{Causley2013a}, we pass to the high-frequency limit.

We obtain the Von Neumann polynomial $\rho^{2}-(2-\beta^{2})\rho+(1-\epsilon)$. We can check that the roots of this polynomial will be complex if $0 < \beta \leq \sqrt{2+2\sqrt{1-\epsilon}}$, and that in this case the roots satisfy

\begin{align}
	|\rho|^{2} &= \frac{1}{4}\left((2-\beta^{2})^{2}+4(1-\epsilon)-(2-\beta^{2})\right) \\
	&= 1-\epsilon < 1
\end{align}

which shows both the stability and dissipative nature of the scheme. $\qed$

We note that the maximum allowed value of $\beta$ is slightly smaller than what is allowed by the corresponding scheme without dissipation. A more detailed analysis shows that the effective damping rate is $\left(\frac{k^{2}}{k^{2}+\alpha^{2}}\right)^{2} \epsilon$, meaning that high frequencies are more rapidly damped than low frequencies.

\subsubsection{Artificial Dissipation in 2D}

Using the notation defined in \cite{Causley2013a}, and again ignoring sources, we have the second order scheme with dissipation,

\begin{align}
	u^{n+1}-2u^{n}+u^{n-1} = -\beta^{2}\mathcal{C}[u^{n}]+ \epsilon \mathcal{C}^{2}[u^{n-1}],
\end{align}

where now $\mathcal{D} = 1-\mathcal{L}_{x}^{-1}\mathcal{L}_{y}^{-1}$ and $\mathcal{C} = \mathcal{L}_{y}^{-1}\mathcal{D}_{x}+\mathcal{L}_{x}^{-1}\mathcal{D}_{y}$. For further details on these operators, see \cite{Causley2013a}. Numerical experiments indicate that the 2D scheme with artificial dissipation is indeed unconditionally stable with the same maximum value of $\beta$ as with the 1D schemes.

\section{Boundary conditions}
\label{sec:BC}

We will now discuss the implementation of boundary conditions. Since our algorithm is dimensionally split, we first develop the boundary conditions in one spatial dimension where the solution \eqref{eqn:2D_explicit} reduces to
\begin{equation}
\label{eqn:Integral_Solution_Full}
u^{n+1}-2u^n+u^{n-1} = -\beta^2 \mathcal{D}_{x}[u^n] + \beta^2\mathcal{L}_{x}^{-1}\left[\frac{1}{\alpha^2} S^n \right](x), \qquad a\leq x \leq b,
\end{equation}
and where we consider the following boundary conditions
\begin{align}
	\label{eqn:Dirichlet}
	\text{Dirichlet:}\qquad & u(a,t) = U_L(t), \quad u(b,t) = U_R(t), \\
	\label{eqn:Neumann}
	\text{Neumann:}\qquad & u_x(a,t) = V_L(t), \quad u_x(b,t) = V_R(t), \\
	\label{eqn:Periodic}
	\text{Periodic}\qquad & u(a,t) = u(b,t), \quad u_x(a,t) = u_x(b,t), \\
	\label{eqn:outflow}
	\text{Outflow:}\qquad & u_t(a,t) =c u_x(a,t), \quad  u_t(b,t)=-c u_x(b,t).
\end{align}
Once these boundary conditions have been derived, we use them to build a boundary solver in 2D.

\subsection{Boundary conditions in one dimension}
\label{sec:bc_1d}
In 1D, this homogeneous solution requires the determination of two coefficients from the imposed boundary conditions and the endpoint values of the particular (integral) solution by solving a $2\times2$ linear system. We now show how to impose several common boundary conditions in 1D. These methods are extended to the 2D case in Section \ref{sec:bc_2d}.

\subsubsection{1D Dirichlet boundary conditions}
Let us begin with Dirichlet boundary conditions \eqref{eqn:Dirichlet}. Evaluating the semi-discrete solution \eqref{eqn:Integral_Solution_Full} at $x = a$ and $b$, we find
\begin{align*}
	U_L(t_{n+1})= 2U_L(t_n)-U_L(t_{n-1}) - \beta^2 \left(U_L(t_n) - I\left[u^n+\frac{1}{\alpha^2}S^n \right](a) - A  - B e^{-\alpha(b-a)}\right), \\
	U_R(t_{n+1})= 2U_R(t_n)-U_R(t_{n-1}) - \beta^2 \left(U_R(t_n) - I\left[u^n+\frac{1}{\alpha^2}S^n \right](b)-A e^{-\alpha(b-a)} - B\right),
\end{align*}
which, after solving for the unknown coefficients can be written as
\begin{align*}
	A^n +\mu	B^n	&= -w_a^D, \\
	\mu	A^n +	B^n	&=  -w_b^D,
\end{align*}
with
\begin{align*}
	w_a^D &= I\left[u^n+\frac{1}{\alpha^2}S^n \right](a) -U_L(t^{n}) - \frac{U_L(t^{n+1}) - 2 U_L(t^{n}) +U_L(t^{n-1}) }{\beta^2}, \\
	w_b^D &= I\left[u^n+\frac{1}{\alpha^2}S^n \right](b) -U_R(t^{n}) -\frac{U_R(t^{n+1}) - 2 U_R(t^{n}) +U_R(t^{n-1}) }{\beta^2},
\end{align*}
and $\mu = e^{-\alpha(b-a)}$. Homogeneous boundary conditions are recovered upon setting $U_L(t) = U_R(t) = 0$. Solving the resulting linear system for the unknowns $A^n$ and $B^n$ gives
\begin{align}
	\label{eqn:wh_Dirichlet}
	A = -\left(\frac{w_a^D - \mu w_b^D }{1-\mu^2}\right), \quad B = -  \left(\frac{w_b^D - \mu w_a^D}{1-\mu^2}\right).
\end{align}

\subsubsection{1D Neumann boundary conditions}
For Neumann conditions, first observe that all dependence on $x$ in the integral solution \eqref{eqn:L_Inverse} is on the Green's function, which is a simple exponential function. Using this, we obtain the following identities
\begin{align}
	\label{eqn:DtN}
	I'(a) = \alpha I(a), \quad  I'(b) = -\alpha I(b).
\end{align}
Now, differentiating the semi-discrete solution \eqref{eqn:Integral_Solution_Full}, and applying the Neumann boundary conditions \eqref{eqn:Neumann} at $x=a$ and $b$ yields
\begin{align*}
	V_L(t_{n+1})=2 V_L(t_n)-V_L(t_{n-1}) - \alpha \beta^2 \left(\frac{1}{\alpha} V_L(t_n) - I\left[u^n+\frac{1}{\alpha^2}S^n \right](a) + A  - B e^{-\alpha(b-a)}\right), \\
	V_R(t_{n+1})=2 V_R(t_n)-V_R(t_{n-1}) - \alpha \beta^2 \left(\frac{1}{\alpha} V_R(t_n) + I\left[u^n+\frac{1}{\alpha^2}S^n \right](b) + A e^{-\alpha(b-a)} - B\right),
\end{align*}
which, after solving for the unknown coefficients can be written as
\begin{align*}
	A^n - \mu	B^n	&= w_a^N, \\
	-\mu	A^n +	B^n	&= w_b^N,
\end{align*}
with
\begin{align*}
	w_a^N &= I\left[u^n+\frac{1}{\alpha^2}S^n \right](a) -\frac{1}{\alpha}V_L(t^n) - \frac{V_L(t^{n+1}) -2 V_L(t^{n}) +V_L(t^{n-1})}{\alpha\beta^2}, \\
	w_b^N &= I\left[u^n+\frac{1}{\alpha^2}S^n \right](b) +\frac{1}{\alpha}V_R(t^n) + \frac{V_R(t^{n+1}) -2 V_R(t^{n}) +V_R(t^{n-1})}{\alpha\beta^2}.
\end{align*}
Upon solving the linear system we obtain
\begin{align}
	\label{eqn:wh_Neumann}
	A = \left(\frac{w_a^N +\mu w_b^N}{1-\mu^2}\right), \quad B =  \left(\frac{w_b^N+\mu w_a^N}{1-\mu^2}\right).
\end{align}

\begin{remark}
	The cases of applying mixed boundary conditions at $x=a$ and $b$ are not considered here, but the details follow from an analogous procedure to that demonstrated above.
\end{remark}
\subsubsection{1D Periodic boundary conditions}
We impose periodic boundary conditions, by assuming that 
\[
u^n(b) = u^n(a), \quad u^n_x(a) = u^n_x(b), \quad n\geq 0.
\]
Enforcing this in the semi-discrete solution \eqref{eqn:Integral_Solution_Full} then yields
\begin{align*}
	I[u^n](a) + A + B\mu  &= I[u^n](b) + A\mu + B, \\
	\alpha \left(I[u^n](a) - A + B\mu\right)  &= \alpha\left(- I[u^n](b) - A\mu + B\right),
\end{align*}
where we have used the identity \eqref{eqn:DtN} applied to derivatives of $I$. Solving this linear system is accomplished quickly by dividing the second equation by $\alpha$, and either adding or subtracting it from the first equation, to produce
\begin{align}
	\label{eqn:wh_Periodic}
	A = \frac{I[u^n](b)}{1-\mu}, \quad B = \frac{I[u^n](a)}{1-\mu}.
\end{align}

\subsubsection{1D Outflow boundary conditions}
\label{sec:Outflow}
When computing wave phenomena, whether we are interested in finite or infinite domains, it is often the case that we must restrict our attention to some smaller subdomain $\Omega$ of the problem, which does not include the physical boundaries. We say that $\Omega$ is the \textit{computational domain}, and that the boundary $\partial \Omega$ is the non-physical, or \textit{artificial boundary}. Under these circumstances, it is necessary to enforce an outflow, or non-reflecting boundary condition, which allows the wave to leave the computational domain, but not incur (non-physical) reflections at the artificial boundary.

For this reason, let us consider the free space solution 
\[
\mathcal{L}^{-1}[u](x) = \frac{\alpha}{2} \int_{-\infty}^\infty u(y)e^{-\alpha|x-y|}dy,
\]
but where we are only interested in evaluating this expression for $x \in \Omega = [a,b]$. Then the contributions can be decomposed as
\begin{align}
	\mathcal{L}^{-1}[u](x)  =& I[u](x) + \frac{\alpha}{2}\int_{-\infty}^a u(y)e^{-\alpha(x-y)}dy + \frac{\alpha}{2} \int_b^\infty u(y) e^{-\alpha(y-x)}dy \nonumber \\
	\label{eqn:Transmission}
	=& I[u](x) + A e^{-\alpha(x-a)} + B e^{-\alpha(b-x)},
\end{align}
where the homogeneous coefficients are
\begin{align}
	A =& \frac{\alpha}{2}\int_{-\infty}^a u(y)e^{-\alpha(a-y)}dy, \\
	B =& \frac{\alpha}{2}\int_b^\infty u(y) e^{-\alpha(y-b)}dy,
\end{align}
which we observe do not depend on $x$. Since the coefficients $A$ an $B$ are the contributions of the integral to the left and right of $[a,b]$ respectively, they can be thought of as transmission conditions (rather than boundary conditions). We make use of this fact to develop outflow boundary conditions, and it will serve as a key idea in our planed follow-up work on a domain decomposition algorithm and multi-core computing with our implicit wave solver.

For the one-dimensional wave equation the exact outflow boundary conditions \eqref{eqn:outflow} turn out to be local in space and time. We emphasize that this is only the case in one spatial dimension, but we shall utilize this fact to obtain an outflow boundary integral solution from the integral equation \eqref{eqn:Transmission}. We extend the support of our function to $(-\infty,\infty)$, and extend the definition of the outflow boundary conditions to the domains exterior to $[a,b]$
\begin{align}
	u_t+cu_x = 0, \quad x\geq b, \\
	u_t-cu_x = 0, \quad x\leq a.
\end{align}

Next, assume the initial conditions have some compact support; for simplicity we will take this support to be $\Omega_0 = [a,b]$. Then after a time $t=t_n$, the domain of dependence of $u^n(x)$ is $\Omega_t = [a-ct_n,b+c t_n]$, since the propagation speed is $c$. Now the free space solution \eqref{eqn:Transmission} becomes
\begin{align}
	\mathcal{L}^{-1}[u^n](x)	=& \frac{\alpha}{2} \int_{a-c t_n}^{b+c t_n} e^{-\alpha|x-y|} u^n(y) dy \nonumber \\
	=& I[u^n](x)+ A^ne^{-\alpha(x-a)} + B^ne^{-\alpha(b-x)}
\end{align}
with coefficients
\begin{align}
	\label{eqn:A_Out_Def}
	A^n &= \frac{\alpha}{2}\int_{a-ct_n}^{a} e^{-\alpha(a-y)} u^n(y) dy, \\
	\label{eqn:B_Out_Def}
	B^n &= \frac{\alpha}{2}\int_{b}^{b+ct_n} e^{-\alpha(y-b)} u^n(y) dy.
\end{align}

At first glance, these coefficients are not at all helpful, as they require computing integrals along spatial domains which not only are outside of the computational domain, but also grow linearly in time. However, we will now make use of the extended boundary conditions to turn these spatial integrals into time integrals, which exist at precisely the endpoints $x=a$ and $b$ respectively. Consider first $x>b$. By assumption, this region contains only right traveling waves, $u(x,t) = u(x-ct)$, and by tracing backward along a characteristic ray we find
\[
u(b+y,t) = u\left(b,t-\frac{y}{c}\right), \quad y>0.
\]
Thus,
\begin{align*}
	B^n &= \frac{\alpha}{2} \int_{0}^{ct_n} e^{-\alpha y} u(b+y,t_n) dy \\
	&= \frac{\alpha c}{2} \int_{0}^{t_n} e^{-\alpha c s} u\left(b,t_n-s\right) ds
\end{align*}
and so $B^n$ is equivalently represented by a convolution in time, rather than space. Now, knowing the history of $u$ at $x=b$ is sufficient to impose outflow boundary conditions. Furthermore, we find in analog to equation \eqref{eqn:IR_def}, a temporal recurrence relation due to the exponential
\begin{align*}
	B^n	&= \frac{\alpha c}{2} \int_{0}^{\Delta t} e^{-\alpha c s} u\left(b,t_n-s\right) ds + e^{-\alpha c \Delta t}\left(\frac{\alpha c}{2}\int_0^{t_{n-1}} e^{-\alpha c s} u\left(b,t_{n-1}-s\right) ds \right) \\
	&= \frac{\beta}{2} \int_{0}^{1} e^{-\beta z} u\left(b,t_n-z\Delta t\right) dz +e^{-\beta} B^{n-1},
\end{align*}
where $\beta = \alpha c \Delta t$, by definition \eqref{eqn:HelmholtzL}. Thus, the coefficient $B^n$, which imposes an outflow boundary condition at $x=b$, can be computed locally in both time and space. To maintain second order accuracy, we fit $u$ with a quadratic interpolant
\[
u(b,t_n-z\Delta t) \approx p(z) = u^n(b) -\frac{z}{2}\left(u^{n+1}(b)-u^{n-1}(b)\right) + \frac{z^2}{2}\left(u^{n+1}(b)-2u^n(b)+u^{n-1}(b)\right)
\]
and integrate the expression analytically 
to arrive at
\begin{equation}
\label{eqn:Outflow_Update_Intermediate}
B^n = e^{-\beta}B^{n-1} +\gamma_0 u^{n+1}(b)+ \gamma_1 u^n(b)+ \gamma_2 u^{n-1}(b)
\end{equation}
where
\begin{align*}
	\gamma_0 =& \frac{E_2(\beta)-E_1(\beta) }{4} = \frac{(1-e^{-\beta})}{2\beta^2}-\frac{(1+e^{-\beta})}{4\beta} \\
	\gamma_1 =& \frac{E_0(\beta)-E_2(\beta) }{2} = -\frac{(1-e^{-\beta})}{\beta^2}+ \frac{1}{\beta}e^{-\beta} +\frac{1}{2}  \\
	\gamma_2 =& \frac{E_2(\beta)+E_1(\beta)}{4} = \frac{(1-e^{-\beta})}{2\beta^2}+\frac{(1-3e^{-\beta})}{4\beta} - \frac{e^{-\beta}}{2}.
\end{align*}
In this outflow update equation \eqref{eqn:Outflow_Update_Intermediate}, the quantities $u^{n+1}(b)$ and $B^n$ are both unknown. In order to determine these values, we must also evaluate the update equation for $u^{n+1}$ \eqref{eqn:Integral_Solution_Full} at $x=b$
\[
u^{n+1}(b) = 2u^{n}(b)-u^{n-1}(b) + \beta^2\left(-u^n(b) + I[u^n](b)+A^n \mu + B^n\right),	\quad \mu = e^{-\alpha(b-a)}.
\]
We now use these two equations to solve for $u^{n+1}(b)$, and eliminate it from the outflow update equation \eqref{eqn:Outflow_Update_Intermediate}, so that
\begin{equation}
\label{eqn:Outflow_Update_B}
-\Gamma_0 \mu A^n + (1-\Gamma_0)B^n = e^{-\beta}B^{n-1} +\Gamma_0 I[u^n](b) +\Gamma_1u^n(b)+ \Gamma_2 u^{n-1}(b)
\end{equation}
where
\[
\Gamma_0 = \beta^2 \gamma_0, \quad
\Gamma_1 = \gamma_1-\gamma_0(\beta^2-2),	\quad
\Gamma_2 = \gamma_2-\gamma_0
\]

\begin{remark}
	While this procedure could be avoided by omitting $u^{n+1}(b)$ in the interpolation stencil, it turns out to be necessary to obtain convergent outflow boundary conditions.
\end{remark}
Likewise, upon considering $x<a$, we find
\begin{equation}
\label{eqn:Outflow_Update_A}
(1-\Gamma_0) A^n -\Gamma_0 \mu B^n = e^{-\beta}A^{n-1} +\Gamma_0 I(a) +\Gamma_1u^n(a)+ \Gamma_2 u^{n-1}(a).
\end{equation}
Solving the resulting linear system produces
\begin{equation}
\label{eqn:Outflow_Update}
A^n = \frac{(1-\Gamma_0)w_a^{\text{Out}} +\mu \Gamma_0w_b^{\text{Out}}} {(1-\Gamma_0)^2-(\mu \Gamma_0)^2}, \quad
B^n = \frac{(1-\Gamma_0)w_b^{\text{Out}} +\mu \Gamma_0w_a^{\text{Out}}} {(1-\Gamma_0)^2-(\mu \Gamma_0)^2},
\end{equation}
where
\begin{align*}
	w_a^{\text{Out}} &= e^{-\beta}A^{n-1} +\Gamma_0 I[u^n](a) +\Gamma_1u^n(a)+ \Gamma_2 u^{n-1}(a), \\
	w_b^{\text{Out}} &= e^{-\beta}B^{n-1} +\Gamma_0 I[u^n](b) +\Gamma_1u^n(b)+ \Gamma_2 u^{n-1}(b)
\end{align*}

\subsection{Boundary conditions in two dimensions}
\label{sec:bc_2d}

We now describe our approach for imposing boundary conditions up to second-order accuracy in two dimensions for our MOL$^T$ formulation of the wave equation.  Boundary conditions must be supplied for the intermediate sweep variable $w$. Since $w = u + O(c^{2}\Delta t^{2})$, for second order accuracy it suffices for $w$ to inherit the boundary condition imposed on the main solution variable $u$. In the case of rectangular, grid-aligned boundaries, the 1D boundary correction terms can be imposed in a line-by-line fashion. Typical applications of periodic and outflow boundary conditions can be imposed in this manner. In the case of complex boundary geometries that are not grid-aligned, Dirichlet boundary conditions can be imposed in a similar line-by-line fashion (by including irregular boundary points as the end points of the sweep lines), but Neumann boundary conditions require more careful treatment due to the resulting coupling of grid lines. Practically speaking,  we only anticipate needing to impose Dirichlet and Neumann boundary conditions on complex boundary geometries, with outflow and periodic boundary conditions being imposed in a line-by-line approach using the 1D results from section \ref{sec:bc_1d} on rectangular domains. Hence, we will limit our discussion here to the implementation of Dirichlet and Neumann boundary conditions in two dimensions. 

\subsubsection{Dirichlet boundary conditions in two dimensions}
\begin{figure}[hbtp]
	\centering
	\subfigure[$x$ lines]{\label{fig:sweep-a}\includegraphics[width=.42\textwidth]{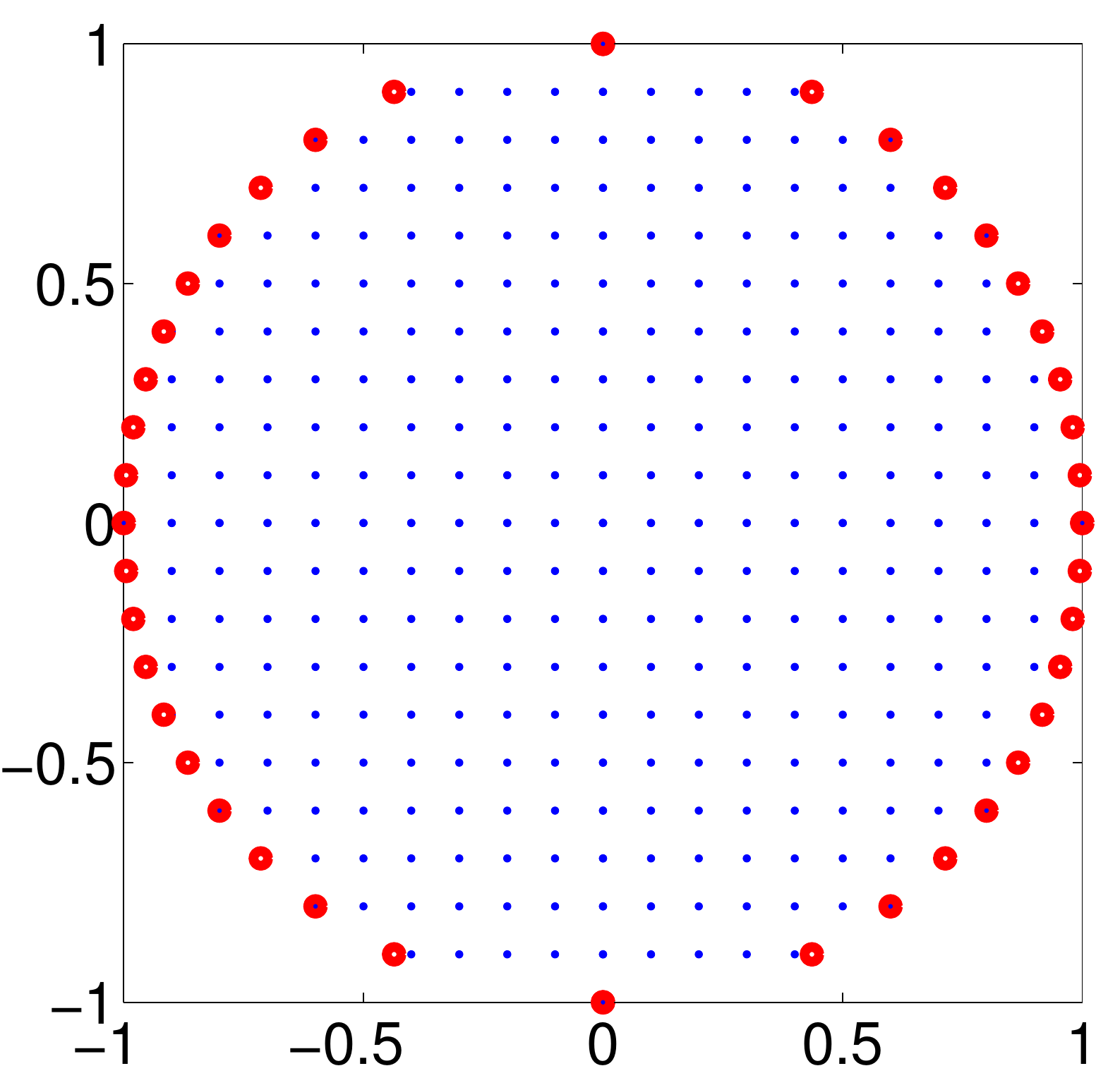}}
	\subfigure[$y$ lines]{\label{fig:sweep-b}\includegraphics[width=.42\textwidth]{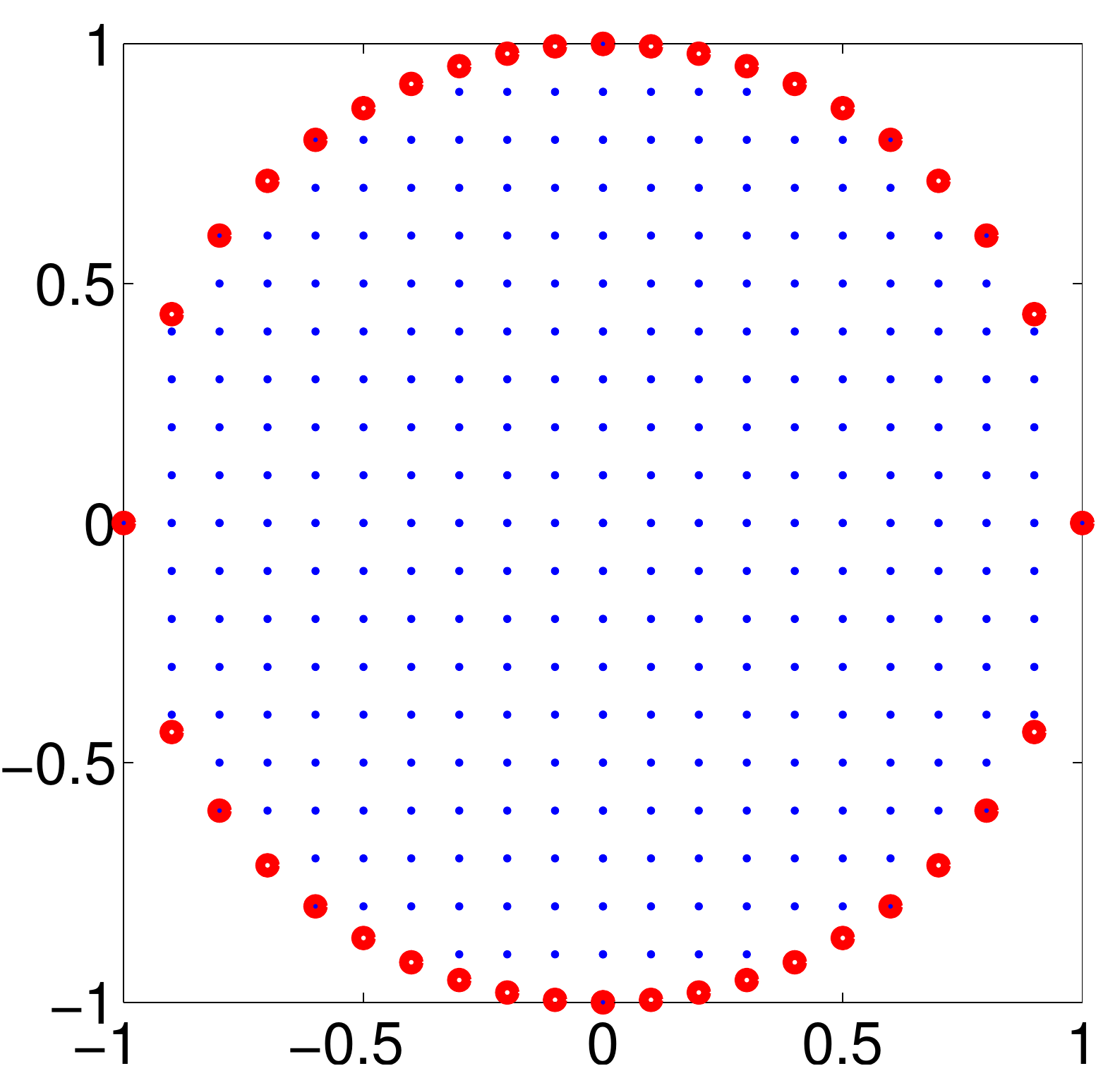}} 
	\caption{Mesh lines used by (a) $x$ integration \eqref{eqn:x_sweep}, and (b) $y$ integration \eqref{eqn:y_sweep} for a circle. The red dots are placed on the boundary, and close each line.}
	\label{fig:sweep}
\end{figure}

Discretization of a general smooth domain $\Omega$ is accomplished by embedding it in a regular Cartesian mesh of say $N_y$ horizontal ($x$) lines and $N_x$ vertical ($y$) lines, and additionally incorporating the termination points of each line, which will lie on the boundary. For example, the lines and boundary points for a circle are shown in Figure \ref{fig:sweep}.

Thus, a line $y = y_k$, which is discretized in $x$, and has endpoints $a_k$ and $b_k$, determined as the intersection of the line with the boundary curve. Analogously we define endpoints of the line $x = x_j$ as $c_j$ and $d_j$, and now have
\begin{align}
	\label{eqn:x_sweep}
	&\mathcal{L}_x^{-1}[u]_k(x)= \frac{\alpha}{2}\int_{a_k}^{b_k} e^{-\alpha|x-x'|} u(x',y)dx' + A_k e^{-\alpha(x-a_k)} + B_k e^{-\alpha(b_k-x)}, \\
	\label{eqn:y_sweep}
	&\mathcal{L}_y^{-1}[u]_j(y) = \frac{\alpha}{2}\int_{c_j}^{d_j} e^{-\alpha|y-y'|} u(x,y')dy' + C_j e^{-\alpha(y-c_j)} + D_j e^{-\alpha(d_j-y)},
\end{align}
for $1\leq k \leq N_y$ and $1\leq j \leq N_x$ respectively. If the boundary is defined using a level set function $C(x,y) = 0$, then the points $x = a_k, b_k$ can be found by solving $C(x,y_k) = 0$, and following from the analogous approach, the endpoints $y = c_j, d_j$ corresponding to $x = x_j$ are found.

Once the domain has been discretized and the lines have been defined, it remains to compute the discrete form of the scheme \eqref{eqn:2D_explicit}. Since the one-dimensional convolution algorithm presented in Section \ref{sec:Fast} is formulated for non-uniform grid points, the embedded boundary points do not affect the implementation of the horizontal \eqref{eqn:x_sweep} and vertical \eqref{eqn:y_sweep} sweeps. Thus, the only point of consideration that remains is that of boundary conditions.

First we consider an intermediate variable $w^{(1)}(x,y)$, defined by
\[
w^{(1)}: = \mathcal{L}_y \left[u^{n+1}-2u^n+u^{n-1}+\beta^2 u^n\right],
\]
and which, according to the implicit scheme \eqref{eqn:Helmholtz_Equation_2D}, can be seen to satisfy
\[
\mathcal{L}_x[w^{(1)}] = \beta^2\left(u^n + \frac{1}{\alpha^2} S^n\right).
\]
Boundary conditions must be applied to $w^{(1)}$ at the boundary points which terminate each line $y = y_k$, defined as $(a_k,y_k)$ and $(b_k,y_k)$. Since $u^{n+1},$ $u^{n},$ and $u^{n-1},$ will be prescribed at these boundary points, we see that the boundary condition will be of the form
\[
w^{(1)}(a_k,y_k) = \lim_{(x,y)\to (a_k,y_k)} \left(1 - \frac{\partial_{yy}}{\alpha^2}\right)\left(u^{n+1}-2u^n+u^{n-1}+\beta^2 u^n\right)(x,y),
\]
and the order of the limit and the partial derivatives commute only when partial derivatives of the boundary data can be constructed. But this is the case only when the tangent line at the boundary is in the $y$-direction, which is precisely why we must restrict our attention to the aforementioned cases for boundary conditions.

Upon introducing a second intermediate variable $w^{(2)}(x,y)$, defined by
\[
w^{(2)}: = u^{n+1}-2u^n+u^{n-1}+\beta^2 u^n,
\]
we make the observation that
\[
\mathcal{L}_y[w^{(2)}] = w^{(1)}.
\]
Boundary conditions are now applied to $w^{(2)}$ along the lines $x = x_j$, at the boundary points $(x_j,c_j)$ and $(x_j,d_j)$, where now no difficulties remain in considering
\[
w^{(2)}(x_j,c_j) = \lim_{(x,y) \to (x_j,c_j)} \left(u^{n+1}-2u^n+u^{n-1}+\beta^2 u^n\right)(x,y),
\]
since the right hand side will be fully prescribed. We summarize this procedure in Algorithm \ref{alg3}.

\begin{algorithm}[htp]
	\caption{Application of Boundary Conditions in 2 dimensions}
	\label{alg3}
	\begin{enumerate}[itemsep=2mm]
		
		\item Initialize temporary variables $w^{(1)}$ and $w^{(2)}$, which are the same size $u^{n}$.
		
		\item For each horizontal line $y = y_k$, for $1\leq k \leq N_y$, create the temporary variable $w^{(1)}_{k}(x)$, defined by
		\[
		\mathcal{L}_x[w^{(1)}_{k}](x) = \beta^2 \left(u^n +\frac{1}{\alpha^2}S^n\right)(x,y_k).
		\]
		The homogeneous coefficients are determined by the fact that
		\[
		w^{(1)} = \mathcal{L}_y[u^{n+1} - 2u^n + u^{n-1} + \beta^2 u^n ],
		\]
		which is applied at $x = a_k, b_k$.
		\item For each vertical line $x = x_j$, for $1\leq j \leq N_x$, create the temporary variable $w^{(2)}_{j}(y)$, defined by
		\[
		\mathcal{L}_x[w^{(2)}_{j}](y) = w^{(1)}_{k}(x).
		\]
		The homogeneous coefficients are determined by the fact that
		\[
		w^{(2)} = [u^{n+1} - 2u^n + u^{n-1} + \beta^2 u^n ],
		\]
		which is applied at $y = c_j, d_j$.
		
		\item Solve for the update
		\[
		u^{n+1} = 2u^n - u^{n-1} -\beta^2 u^n +w^{(2)}.
		\]
		
		\item The dimensional splitting error is corrected by adding
		\[
		\beta^2\mathcal{D}_x\mathcal{D}_y[u^n] = \beta^2\mathcal{D}_y[u^n] + w^{(1)}-w^{(2)}.
		\]

	\end{enumerate}
\end{algorithm}

\subsubsection{Neumann boundary conditions in two dimensions}
\label{sec:embedded_neumann_bc}

In implementing Neumann boundary conditions for boundary geometries conforming to grid lines, such as a rectangular domain, we can directly impose a two-point boundary correction. One way to extend this method to a general polygonal domain would be to use multiple overset grids, each aligned with a boundary segment, which communicate with the interior grid through interpolation on a ghost cell region, though we do not pursue that approach in this work.

For curved boundaries, an alternative approach is that of an embedded boundary method, which involves determining the Dirichlet values at the endpoints of each $x$- and $y$-sweep lines that result in the approximate satisfaction of the Neumann boundary condition (in effect, constructing an approximate Neumann-to-Dirichlet map). We present the implementation of an embedded boundary method for Neumann boundary conditions for the implicit wave solver on a curved boundary geometry. The approach taken here follows the work in \cite{kreiss2004difference}, which proposes an embedded boundary method for Neumann boundary conditions with a finite difference method for the wave equation. The analysis in that work suggests that, on the continuous level, the modified equations and boundary conditions resulting from typical truncation error terms possess unstable boundary layer solutions, so that the addition of a dissipative term is necessary to achieve a stable method. This is consistent with our experience in the implementation described here, with the embedded boundary method becoming unstable when applied to the non-dissipative dispersive solver, but remaining stable for the diffusive solver, which is dissipative. 

In the following, we briefly describe the two-point boundary correction method for a grid-aligned rectangular boundary. We then describe the embedded boundary method for a 1D problem. This method requires an iterative procedure, which we show, in the setting of the 1D problem, to be a convergent contraction mapping with a rate of convergence that depends on the CFL number. We describe the implementation of the embedded boundary method in 2D, and finally give numerical results.

\subsubsection{Description of the Two-Point Boundary Correction Method}

In a rectangular domain where the boundaries conform to grid lines, it is straightforward to impose the two-point boundary correction terms in a line-by-line fashion, since in this case, the grid lines are not coupled through the normal derivative. As this is a simple extension of the 1D boundary correction algorithm, we do not elaborate further.

\subsubsection{Description of the Embedded Boundary Method and Proof of Convergence of the Iterative Solution in 1D}

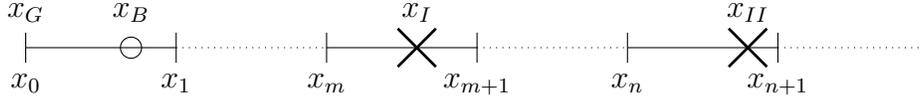
\begin{figure}[hbh] 
	\begin{center}
		\begin{tikzpicture}[scale = 2]
		
		\draw (-1.5,0)--(-.5,0);
		\draw[dotted] (-.5,0)--(.5,0);
		\draw (.5,0)--(1.5,0);
		\draw[dotted] (1.5,0)--(2.5,0);
		\draw (2.5,0)--(3.5,0);
		\draw[dotted] (3.5,0)--(4.5,0);
		
		\draw (-1.5,-.1)--(-1.5,.1);
		\draw (-.5,-.1)--(-.5,.1);
		\draw (.5,-.1)--(.5,.1);
		\draw (1.5,-.1)--(1.5,.1);
		\draw (2.5,-.1)--(2.5,.1);
		\draw (3.5,-.1)--(3.5,.1);
		
		\node (c) at (1.1,0) {\Cross};
		\node (c) at (3.3,0) {\Cross};
		
		
		
		
		\draw (-1.5,-.1) node[anchor=north] {$x_{0}$};
		\draw (-1.5,.1) node[anchor=south] {$x_{G}$};
		\draw (-0.8,.1) node[anchor=south] {$x_{B}$};
		\draw (-0.8,0) circle [radius = 2pt];
		\draw (-.5,-.1) node[anchor=north] {$x_{1}$};
		\draw (.5,-.1) node[anchor=north] {$x_{m}$};
		\draw (1.5,-.1) node[anchor=north] {$x_{m+1}$};
		
		\draw (1.1,.1) node[anchor=south] {$x_{I}$};
		
		\draw (2.5,-.1) node[anchor=north] {$x_{n}$};
		\draw (3.5,-.1) node[anchor=north] {$x_{n+1}$};
		\draw (3.3,.1) node[anchor=south] {$x_{II}$};

		\end{tikzpicture}
	\end{center} 
	\caption{Boundary geometry in 1D.}
	\label{fig:boundary_geom_1d}
\end{figure}

We consider the situation of a one-dimensional domain $\left\lbrace x_{B}<x \right\rbrace$ with a single boundary point not aligned with the grid points, as displayed in Figure \ref{fig:boundary_geom_1d}. We have grid points $x_{0},x_{1},...$ with uniform grid spacing $x_{i+1}-x_{i} = \Delta x$, boundary location $x_{B}$, and ghost point location $x_{G}=x_{0}$. We define interior points to be any grid points lying within the domain (including the boundary), and exterior points to be any grid points lying outside of the domain. We define a ghost point to be any exterior point for which at least one of the neighboring points $x_{i \pm 1}$ is an interior point.  We neglect the right boundary in the present analysis for simplicity, though it can be extended to the case with both boundaries. We consider applying the diffusive version of the wave solver, having calculated the convolution integral $I(x)$, and now needing to find the value of the coefficient $A$ such that the solution $u(x)$ at the next time step is given by 
\begin{align*}
	u(x) = I(x) + A e^{-\alpha (x-x_{G})}, \quad x \geq x_{0}
\end{align*}

Given the value of the convolution integral and the solution at the ghost point, $I_{G} = I(x_{G})$ and $u_{G}=u(x_{G})$, respectively, the coefficient may be computed as $A =  u_{G}-I_{G}$. We now describe the procedure for determining the value of the solution at the ghost point, $u_{G}$, that leads to a solution consistent with homogeneous Neumann boundary conditions to second-order accuracy. We construct a quadratic interpolant using the boundary condition and interior interpolation points $x_{I} = x_{B}+ \Delta s_{I}$ and $x_{II} = x_{B}+2\Delta s_{I}$, lying in between grid points $x_{m}$ and $x_{m+1}$, and $x_{n}$ and $x_{n+1}$, respectively. The interpolation distance will be chosen such that $\Delta x < \Delta s_{I} < (3/2) \Delta x$. We define the distances $\xi_{G} = x_{B}-x_{G}$, $\xi_{I} = x_{I}-x_{B} = \Delta s_{I}$, and $\xi_{II} = x_{II}-x_{B} = 2\Delta s_{I}$, and construct a quadratic Hermite-Birkhoff \cite{birkhoff1906general} interpolant $P(\xi)$ by imposing the conditions $P'(0) = 0$, $P(\xi_{I}) = u_{I}$ and $P(\xi_{II}) = u_{II}$. We then obtain the following second-order approximation to the ghost point value, given by

\begin{align*}
	u_{G} &= P(\xi_{G}) + O(\Delta x^{2}) = u_{I} \frac{\xi_{II}^{2}-\xi_{G}^{2}}{\xi_{II}^{2}-\xi_{I}^{2}}+ u_{II} \frac{\xi_{G}^{2}-\xi_{I}^{2}}{\xi_{II}^{2}-\xi_{I}^{2}} + O(\Delta x^{2}) \\
	&= \gamma_{I} u_{I}+\gamma_{II} u_{II}.
\end{align*}

As the coefficients $\gamma_{I} = \frac{\xi_{II}^{2}-\xi_{G}^{2}}{\xi_{II}^{2}-\xi_{I}^{2}} > 0$ and $\gamma_{II} = \frac{\xi_{G}^{2}-\xi_{I}^{2}}{\xi_{II}^{2}-\xi_{I}^{2}} < 0$  are $O(1)$, we only need supply second-order accurate approximations to $u_{I}$ and $u_{II}$ to maintain overall second-order accuracy. (The coefficients would be $O(1/\Delta x)$ in the case of nonhomogeneous Neumann boundary conditions, which would require third-order accurate approximations to $u_{I}$ and $u_{II}$. For simplicity, we consider only the case of homogeneous Neumann boundary conditions in the present work.) Such approximations may be obtained through linear interpolation, giving

\begin{align*}
	u_{I} &= \sigma_{I} u_{m} + (1-\sigma_{I}) u_{m+1} \\
	u_{II} &= \sigma_{II} u_{n} + (1-\sigma_{II}) u_{n+1}
\end{align*}

where $\sigma_{I} = \frac{x_{m+1}-x_{I}}{\Delta x}$, $\sigma_{II} = \frac{x_{n+1}-x_{II}}{\Delta x}$, and $u_{j} = u(x_{j})$ are the values of the function at the uniform gridpoints for $j=m,m+1,n,n+1$.

Hence, to determine the ghost point value $u_{G}$ that leads to a solution consistent with homogeneous Neumann boundary conditions, we should solve the following system of equations.

\begin{align*}
	u_{G} &= \gamma_{I} (\sigma_{I} u_{m} + (1-\sigma_{I}) u_{m+1}) + \gamma_{II} (\sigma_{II} u_{n} + (1-\sigma_{II}) u_{n+1}) \\
	u_{m} &= I_{m} + (u_{G}-I_{G})e^{-\alpha(x_{m}-x_{0})} \\
	u_{m+1} &= I_{m+1} + (u_{G}-I_{G})e^{-\alpha(x_{m+1}-x_{0})} \\
	u_{n} &= I_{n} + (u_{G}-I_{G})e^{-\alpha(x_{n}-x_{0})} \\
	u_{n+1} &= I_{n+1} + (u_{G}-I_{G})e^{-\alpha(x_{n+1}-x_{0})} 
\end{align*}

with $\gamma_{I}$, $\gamma_{II}$, $\sigma_{I}$ and $\sigma_{II}$ defined as above, and where $I_{j} = I(x_{j})$ is the convolution integral evaluated at uniform grid points for $j=m,m+1,n,n+1$. This system can be solved formally for $u_{G}$, giving

\begin{align*}
	u_{G} &= \left[ \gamma_{I} \left(\sigma_{I}(I_{m}-I_{G} e^{-\alpha(x_{m}-x_{G})})+(1-\sigma_{I})(I_{m+1}-I_{G} e^{-\alpha(x_{m+1}-x_{G})})\right) + \right. \\
	& \quad \left. \gamma_{II} \left(\sigma_{II}(I_{n}-I_{G} e^{-\alpha(x_{n}-x_{G})})+(1-\sigma_{II})(I_{n+1}-I_{G} e^{-\alpha(x_{n+1}-x_{G})})\right) \right] \div \\
	& \quad\left[ 1- \gamma_{I} \left(\sigma_{I}e^{-\alpha(x_{m}-x_{G})}+(1-\sigma_{I})e^{-\alpha(x_{m+1}-x_{G})}\right) - \right. \\
	& \quad \left. \gamma_{II} \left(\sigma_{II}e^{-\alpha(x_{n}-x_{G})}+(1-\sigma_{II})e^{-\alpha(x_{n+1}-x_{G})}\right) \right]
\end{align*}

To show that this solution formula is well-defined, we argue that 
\begin{align*}
	0< K &:= \gamma_{I} \left(\sigma_{I}e^{-\alpha(x_{m}-x_{G})}+(1-\sigma_{I})e^{-\alpha(x_{m+1}-x_{G})}\right) + \\
	& \quad + \gamma_{II} \left(\sigma_{II}e^{-\alpha(x_{n}-x_{G})}+(1-\sigma_{II})e^{-\alpha(x_{n+1}-x_{G})}\right) < 1
\end{align*}

for the relevant values of $m$, $n$ and $\xi_{G}$, $\xi_{I}$, and $\xi_{II}$. We define $d = e^{-\alpha \Delta x}$, and noting that $0 < d < 1$, $m<n$, $\xi_{G} < \xi_{I} = \Delta s_{I} < \xi_{II} = 2\Delta s_{I}$, $\gamma_{I}>0$, $\gamma_{II}<0$, $0 \leq \sigma_{I} \leq 1$, and $0 \leq \sigma_{II} \leq 1 $, we obtain

\begin{align*}
	K &= \gamma_{I} d^{m} \left[\sigma_{I}+(1-\sigma_{I})d\right] + \gamma_{II} d^{n} \left[\sigma_{II}+(1-\sigma_{II})d\right] \\
	&\leq \gamma_{I} d^{m} + \gamma_{II} d^{n+1} \\
	&= \frac{d^{m} \xi_{II}^{2}-d^{n+1}\xi_{I}^{2}+\xi_{G}^{2}(d^{n+1}-d^{m})}{\xi_{II}^{2}-\xi_{I}^{2}} \\
	&\leq \frac{d^{m} \xi_{II}^{2}-d^{n+1}\xi_{I}^{2}}{\xi_{II}^{2}-\xi_{I}^{2}} \\
	&= \frac{4\Delta s_{I}^{2} d^{m}-d^{n+1} \Delta s_{I}^{2}}{4\Delta s_{I}^{2}-\Delta s_{I}^{2}} = \frac{4d^{m}-d^{n+1}}{3}
\end{align*}

Now, since $\Delta x < \Delta s_{I} < (3/2) \Delta x$, we can see that it is the case that either $m=1$ and $n=2$ or $3$, or that $m=2$ and $n=3$. It is then a matter of some simple calculus to check that that the functions $f_{m,n} (x) = (4x^{m}-x^{n+1})/3$ satisfy $f_{m,n}(x)<1$ for $0<x<1$ and the given combinations of $m$ and $n$. This proves that $K<1$ for the relevant values of the parameters, so that the solution for $u_{G}$ is well-defined. We note, however, that $K$ approaches 1 as $d$ approaches 1, that is, as the CFL number becomes large. Thus, we may expect an ill-conditioned system when the CFL number is very large.

To obtain the lower bound on $K$, we observe
\begin{align*}
	K &= \gamma_{I} d^{m} \left[\sigma_{I}+(1-\sigma_{I})d\right] + \gamma_{II} d^{n} \left[\sigma_{II}+(1-\sigma_{II})d\right] \\
	&\geq \gamma_{I} d^{m+1} + \gamma_{II} d^{n} \\
	&= \frac{d^{m+1}(\xi_{II}^{2}-\xi_{G}^{2})+d^{n}(\xi_{G}^{2}-\xi_{I}^{2})}{\xi_{II}^{2}-\xi_{I}^{2}} \\
	&= \frac{d^{m+1}(4\xi_{I}^{2}-\xi_{G}^{2})+d^{n}(\xi_{G}^{2}-\xi_{I}^{2})}{\xi_{II}^{2}-\xi_{I}^{2}} \\
	&= \frac{(d^{m+1}-d^{n})(\xi_{I}^{2}-\xi_{G}^{2})+3d^{m+1}\xi_{I}^{2}}{\xi_{II}^{2}-\xi_{I}^{2}} > 0 
\end{align*}

In the two-dimensional case, the line-by-line solution method couples the ghost point values, and a general explicit solution formula is impossible to write down. In principle, one may write out and directly solve a linear system to obtain the ghost point values. Instead, we propose an iterative solution method that avoids the formation of a matrix. We now describe this iterative solution method and prove its convergence in the context of the one-dimensional problem described above.

Suppose we have the convolution integral evaluated at the gridpoints, $I_{j}$, and a $k$-th iterate for the ghost point value, $u_{G}^{k}$. Then we may obtain the next iterate by the formulas
\begin{align*}
	u_{m}^{k+1} &= I_{m} + (u_{G}^{k}-I_{G})e^{-\alpha(x_{m}-x_{0})} \\
	u_{m+1}^{k+1} &= I_{m+1} + (u_{G}^{k}-I_{G})e^{-\alpha(x_{m+1}-x_{0})} \\
	u_{n}^{k+1} &= I_{n} + (u_{G}^{k}-I_{G})e^{-\alpha(x_{n}-x_{0})} \\
	u_{n+1}^{k+1} &= I_{n+1} + (u_{G}^{k}-I_{G})e^{-\alpha(x_{n+1}-x_{0})} \\
	u_{G}^{k+1} &= \gamma_{I} (\sigma_{I} u_{m}^{k+1} + (1-\sigma_{I}) u_{m+1}^{k+1}) + \gamma_{II} (\sigma_{II} u_{n}^{k+1} + (1-\sigma_{II}) u_{n+1}^{k+1})
\end{align*}
where quantities are defined as above. Now, to prove the convergence of the interation, we show it is contractive. Taking the difference of two iterates, we have
\begin{align*}
	|u_{G}^{k+1}-u_{G}^{k}| &= |\gamma_{I} (\sigma_{I} (u_{m}^{k+1}-u_{m}^{k}) + (1-\sigma_{I}) (u_{m+1}^{k+1}-u_{m+1}^{k})) + \\
	& \quad \gamma_{II} (\sigma_{II} (u_{n}^{k+1}-u_{n}^{k}) + (1-\sigma_{II}) (u_{n+1}^{k+1}-u_{n+1}^{k}))| \\
	& = |\gamma_{I} \left(\sigma_{I}(u_{G}^{k}-u_{G}^{k-1})e^{-\alpha(x_{m}-x_{G})}+(1-\sigma_{I})(u_{G}^{k}-u_{G}^{k-1})e^{-\alpha(x_{m+1}-x_{G})}\right) + \\
	& \quad + \gamma_{II} \left(\sigma_{II}(u_{G}^{k}-u_{G}^{k-1})e^{-\alpha(x_{n}-x_{G})}+(1-\sigma_{II})(u_{G}^{k}-u_{G}^{k-1})e^{-\alpha(x_{n+1}-x_{G})}\right)| \\
	&\leq K |u_{G}^{k}-u_{G}^{k-1}|
\end{align*}
where $0 < K < 1$ as defined above. Hence, the Contraction Mapping Theorem implies that the iteration converges to a unique fixed point (which is the solution given in the formula above). We note again that $K$ approaches 1 as the CFL number becomes large, so that the rate of convergence will become slower for larger CFL numbers.

\subsubsection{Description of the Method in 2D}

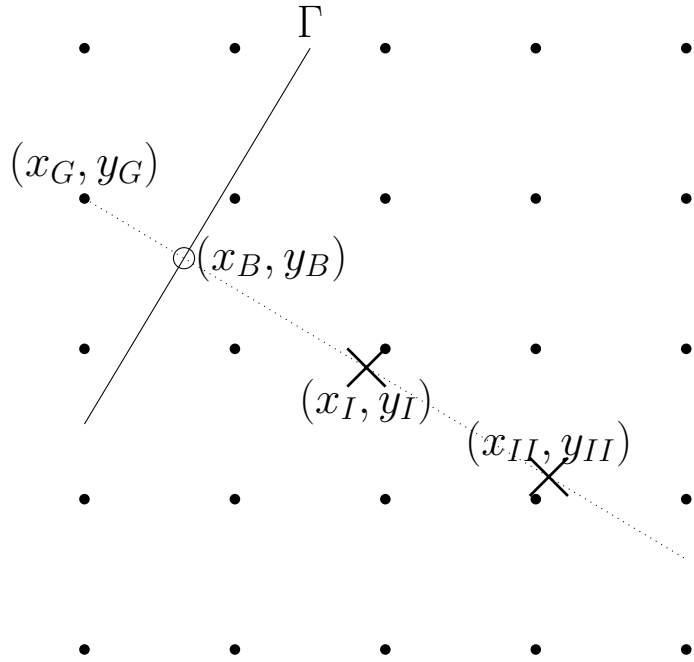
\begin{figure}[hbh] 
	\begin{center}
		\begin{tikzpicture}[scale = 2]
		
		\foreach \x in {0,1,2,3,4}
		\foreach \y in {0,1,2,3,4}
		\draw (\x,\y)  node[fill,circle,inner sep=0pt,minimum size=4pt] {} ;
		
		\draw (0,1.5)--(1.5,4);
		
		\draw (1.5,4) node[anchor=south] {\Large $\Gamma$};
		
		\coordinate (boundaryIntersection) at ({45/68}, {(-3/5)*(45/68)+3});
		
		\coordinate (normalEndpoint) at (4,{(-3/5)*4+3});
		
		\coordinate (intPointI) at ({(5/68)*(9+4*sqrt(17))},{(-3/5)*(5/68)*(9+4*sqrt(17))+3});
		\coordinate (intPointII) at ({(5/68)*(9+8*sqrt(17))},{(-3/5)*(5/68)*(9+8*sqrt(17))+3});
		%
		
		\draw[dotted] (0,3)--(normalEndpoint);
		
		\draw (boundaryIntersection) circle [radius = 2pt];
		
		\node (c) at (intPointI) {\Cross};
		\node (c) at (intPointII) {\Cross};	
		
		\draw (0,3) node[anchor=south] {\Large $(x_{G},y_{G})$};
		\draw (intPointI) node[anchor=north] {\Large $(x_{I},y_{I})$};
		\draw (intPointII) node[anchor=south] {\Large $(x_{II},y_{II})$};
		\draw (boundaryIntersection) node[anchor=west] {\Large $(x_{B},y_{B})$};
		
		\end{tikzpicture}
	\end{center} 
	\caption{Boundary geometry in 2D.}
	\label{fig:boundary_geom_2d}
\end{figure}

We now describe the implementation of the embedded Neumann boundary condition in the 2D case. We consider the situation displayed in Figure \ref{fig:boundary_geom_2d}, in which we need to determine the value of our unknown $u_{G}$ at the ghost point location $(x_{G},y_{G})$. In the 2D case, we define a ghost point to be any exterior point $(x_{i},y_{j})$ for which at least one of the neighboring points $(x_{i \pm 1},y_{j})$ or $(x_{i},y_{j\pm 1})$ is an interior point.  Similarly to the 1D case, we will construct a quadratic Hermite-Birkhoff boundary interpolant $P(\xi)$ along the direction normal to the boundary, which intersects the boundary curve $\Gamma$ at location $(x_{B},y_{B})$, and supply the interior interpolation point values $u_{I}$ and $u_{II}$, at points $(x_{I},y_{I})$ and $(x_{II},y_{II})$, respectively, by further interpolation from interior grid points. These points are selected along the normal, in analogy to the 1D case, such that $\xi_{I} = |(x_{I},y_{I})-(x_{B},y_{B})| = \Delta s_{I}$ and $\xi_{II} = |(x_{II},y_{II})-(x_{B},y_{B})| = 2 \Delta s_{I}$, where we will typically take $\Delta s_{I} = \sqrt{2} \Delta x$. We construct a quadratic Hermite-Birkhoff interpolant $P(\xi)$ by imposing the conditions $P'(0) = 0$, $P(\xi_{I}) = u_{I}$ and $P(\xi_{II}) = u_{II}$. Defining further the distance from the boundary to the ghost point $\xi_{G} = |(x_{G},y_{G})-(x_{B},y_{B})|$, we obtain, as in the 1D case, the following second-order approximation to the ghost point value, given by

\begin{align*}
	u_{G} &= P(\xi_{G}) + O(\Delta x^{2}) = u_{I} \frac{\xi_{II}^{2}-\xi_{G}^{2}}{\xi_{II}^{2}-\xi_{I}^{2}}+ u_{II} \frac{\xi_{G}^{2}-\xi_{I}^{2}}{\xi_{II}^{2}-\xi_{I}^{2}} + O(\Delta x^{2}) \\
	&= \gamma_{I} u_{I}+\gamma_{II} u_{II}.
\end{align*}

where the coefficients are defined as $\gamma_{I} = \frac{\xi_{II}^{2}-\xi_{G}^{2}}{\xi_{II}^{2}-\xi_{I}^{2}} > 0$ and $\gamma_{II} = \frac{\xi_{G}^{2}-\xi_{I}^{2}}{\xi_{II}^{2}-\xi_{I}^{2}} < 0$. In the 2D case, we find approximations to $u_{I}$ and $u_{II}$ through bilinear interpolation. This is in contrast to \cite{kreiss2004difference}, who find the intersection of the normal with grid lines, then interpolate along the grid lines. We have also implemented a second-order accurate version of this approach and compared to the bilinear interpolation scheme proposed here. We have found that the two schemes behave similarly, however the bilinear interpolation scheme is slightly more accurate and simpler to code, not requiring to handle separate cases of intersection with horizontal and vertical grid lines. The bilinear interpolation scheme is standard, but we give it here for completeness. If the interpolation point $u_{I}$ lies in a cell with corners $(x_{i},y_{j})$, $(x_{i+1},y_{j})$, $(x_{i+1},y_{j+1})$ and $(x_{i},y_{j+1})$, then we have the following approximation for $u_{I}$:

\begin{align*}
	u_{I} &= w_{1}u_{i,j}+w_{2}u_{i+1,j}+w_{3}u_{i+1,j+1}+w_{4}u_{i,j+1}
\end{align*}

where $w_{1} = \frac{(x_{i+1}-x_{I})(y_{j+1}-y_{I})}{\Delta x \Delta y}$, $w_{2} = \frac{(x_{I}-x_{i})(y_{j+1}-y_{I})}{\Delta x \Delta y}$, $w_{3} = \frac{(x_{I}-x_{i})(y_{I}-y_{j})}{\Delta x \Delta y}$ and $w_{4} = \frac{(x_{i+1}-x_{I})(y_{I}-y_{j})}{\Delta x \Delta y}$. With this interpolation scheme established, we now outline the algorithm for the 2D dimensionally-split wave solver.

The above interpolation procedure applies regardless of the variety of the wave solver that it is used with, provided that the wave solver has sufficient dissipation to maintain stability. We now describe the rest of the embedded boundary algorithm in the context of the diffusive wave solver, though it it may be similarly applied to the the dispersive scheme with artificial dissipation described above. In analogy to the iteration presented in the 1D case, the 2D iterative algorithm proceeds by using the interpolation scheme to provide values at the ghost points, which in turn provide new values for the boundary correction coefficients, which are then used to update the values at the interior grid points, comprising one full iteration. It should be noted that not all interior grid points need be updated in the iteration, only those near the boundary that lie within the boundary interpolation stencils. 

Using values from previous time steps, the initial guess for the interior grid points in the boundary interpolation stencils may be given by extrapolation in time, as $u^{n+1,0} = 2u^{n}-u^{n-1}$ (linear extrapolation) or $u^{n+1,0} = 3u^{n}-3u^{n-1}+u^{n-2}$ (quadratic extrapolation). Either extrapolated initial guess provides a modest reduction in the number of iterations required versus a zero initial guess, with only a slight further reduction in the number of iterations going from linear to quadratic extrapolation. An effective stopping criterion for iteration is $|u^{n+1,l+1}-u^{n+1,l}|_{\infty} < \epsilon$, where $\epsilon$ is some chosen tolerance, which may be chosen to be quite small, as the iteration is a fixed point interation. In the numerical example, we choose a tolerance of $10^{-15}$, and we achieve convergence in less than 40 iterations at a CFL number of 2.

In applying the diffusive version of the wave solver, we assume we have previous time steps $u^{n}$, $u^{n-1}$ and $u^{n-2}$. We have to solve the modified Helmholtz equation with homogeneous Neumann boundary conditions,

\begin{align*}
	\left(1-\frac{1}{\alpha^{2}}\nabla^{2}\right)u^{n+1} = \frac{1}{2}\left(5u^{n}-4u^{n-1}+u^{n-2}\right) & \mbox{ in } \Omega \\
	\frac{\partial u}{\partial n} = 0 & \mbox{ on } \Gamma = \partial \Omega
\end{align*}

where $\alpha = \frac{\sqrt{2}}{c \Delta t}$. We apply dimensional splitting to find

\begin{align*}
	\left(1-\frac{1}{\alpha^{2}}\nabla^{2}\right)u^{n+1} = \left(1-\frac{1}{\alpha^{2}}\partial_{xx}\right)\left(1-\frac{1}{\alpha^{2}}\partial_{yy}\right)u^{n+1} + O\left(\frac{1}{\alpha^{4}}\right)
\end{align*}

so we define $w  = \left(1-\frac{1}{\alpha^{2}}\partial_{yy}\right)u$, and noting that $w = u+O\left(\left(c\Delta t \right)^{2}\right)$ so that $\frac{\partial w}{\partial n} = \frac{\partial u}{\partial n}+O\left(\left(c\Delta t \right)^{2}\right)$, we obtain the following approximate system

\begin{align*}
	\left(1-\frac{1}{\alpha^{2}}\partial_{xx}\right)w^{n+1} = \frac{1}{2}\left(5u^{n}-4u^{n-1}+u^{n-2}\right) & \mbox{ in } \Omega \\
	\frac{\partial w}{\partial n}^{n+1} = 0 & \mbox{ on } \Gamma = \partial \Omega \\
	\left(1-\frac{1}{\alpha^{2}}\partial_{yy}\right)u^{n+1} = w^{n+1} & \mbox{ in } \Omega \\
	\frac{\partial u}{\partial n}^{n+1} = 0 & \mbox{ on } \Gamma = \partial \Omega
\end{align*}

We now suppose our domain is embedded in a uniform Cartesian grid, with horizontal grid lines corresponding to $y=y_{k}$, $1\leq k \leq N_{y}$ and vertical grid lines corresponding to $x=x_{j}$, $1\leq j \leq N_{x}$. The embedded boundary algorithm will be applied when calculating the intermediate variable $w^{n+1}$ in horizontal line sweeps as well as the solution variable $u^{n+1}$ in vertical line sweeps. The iterations for these two variables are separate; first, the iterative procedure is applied to $w$ to convergence, and then this value of $w$ is used to compute $u$, and the iterative procedure is applied to $u$ to convergence. However, in each iteration, the grid lines are coupled through the interpolation scheme, so that all grid lines must be iterated together. The overall iterative algorithm is described in \ref{alg:eb_algorithm}, with details specified for the iteration on $w$. The iteration on $u$ is very similar, and so we omit the details.

\begin{algorithm}[htp]
	\caption{Application of Neumann Boundary Conditions in 2 dimensions}
	\label{alg:eb_algorithm}
	\begin{enumerate}[itemsep=2mm]
		
		\item \textbf{(Initialization of ghost points)} Perform the interpolation scheme described above to obtain the values of $u^{n}$, $u^{n-1}$, and $u^{n-2}$ at the ghost points, which are the endpoints of the horizontal and vertical grid lines.
		
		\item \textbf{(Evaluation of particular solution)} For each horizontal line $y = y_{k}$, for $1\leq k \leq N_{y}$, with ghost (end) points $x=a_{k}$ and $b_{k}$ find the particular solution $w^{n+1}_{p,k}$ for the intermediate variable $w^{n+1}_{k}(x)$ by evaluating the discrete convolution operator
		\[
		w^{n+1}_{p,k}(x_{j}) =  \frac{\alpha}{4}\int_{a_{k}}^{b_{k}}[5u^{n}-4u^{n-1}+u^{n-2}](x',y_{k})e^{-\alpha|x_{j}-x'|} \, dx'
		\]
		
		for each grid point $x_{j}$ in the horizontal line, including the ghost points.
		
		\item \textbf{(Boundary correction initialization)} For each horizontal line $y = y_{k}$, set the initial guess for the intermediate variable via extrapolation, $w^{n+1,0}_{k} = 3w^{n}_{k}-3w^{n-1}_{k}+w^{n-2}_{k}$, on the interior points within the boundary interpolation stencil.
		
		\item \textbf{(Boundary correction iteration)} For each horizontal line $y = y_{k}$, perform the interpolation scheme using the interior values of $w^{n+1,l}_{p,k}$ to find the ghost point values. Using these ghost point values, apply the boundary correction on each line to obtain the updated intermediate variable,
		
		\[
		w^{n+1,l+1}_{k}(x_{j}) = w^{n+1}_{p,k}(x_{j})+A_{k}e^{-\alpha(x_{j}-a_{k})} + B_{k}e^{-\alpha(b_{k}-x_{j})}
		\]
		
		for the values of $x_{j}$ lying within the boundary interpolation stencil, where $A_{k} = \frac{w^{n+1,l}_{k}(a_{k})-w^{n+1}_{p,k}(a_{k})-\mu_{k} \left(w^{n+1,l}_{k}(b_{k})-w^{n+1}_{p,k}(b_{k})\right)}{1-\mu_{k}^{2}}$, $B_{k} = \frac{w^{n+1,l}_{k}(b_{k})-w^{n+1}_{p,k}(b_{k})-\mu_{k} \left(w^{n+1,l}_{k}(a_{k})-w^{n+1}_{p,k}(a_{k})\right)}{1-\mu_{k}^{2}}$, $\mu_{k} = e^{-\alpha(b_{k}-a_{k})}$. Check for convergence, and if converged, store the intermediate variable $w^{n+1}$.
		
		\item Repeat this process for the vertical line sweeps, using the intermediate variable $w^{n+1}$ to calculate the particular solution for $u^{n+1}$, then apply the bounday correction interation.

	\end{enumerate}
\end{algorithm}

\section{Numerical Results}
\label{sec:results}

\subsection{Double Circle Cavity}

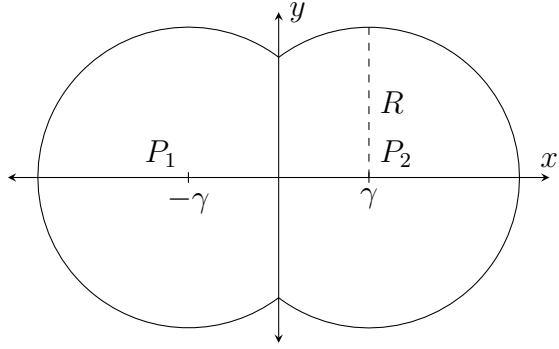
\begin{figure}[H]
	\centering
	\begin{tikzpicture}[scale=2,>=stealth]
	
	\draw[<->] (-1.8,0) -- (1.8,0)
	node [above] {$x$};
	\draw[<->] (0,-1.1) -- (0,1.1)
	node [right] {$y$};
	\draw (0,-.8) arc (-126.9:126.9:1);
	\draw (0,.8) arc (53.1:306.9:1);
	\draw[dashed] (.6,0) -- (.6,1);
	\draw (.6,.5) node[anchor=west] {$R$};
	
	\foreach \x in {-.6,.6}
	{
		\draw (\x,-1pt)--(\x,1pt);
	}
	
	\draw (.6,0) node[anchor=north] {$\gamma$};
	\draw (-.6,0) node[anchor=north] {$-\gamma$};
	
	\draw (.6,0) node[anchor=south west] {$P_{2}$};
	\draw (-.6,0) node[anchor=south east] {$P_{1}$};
	
	\end{tikzpicture}
	\caption{Double circle geometry.}
	\label{fig:dblcirc_geom}
\end{figure}

In this example, we solve the wave equation with homogeneous Dirichlet boundary conditions on a 2D domain $\Omega$ which is, as in Figure \ref{fig:dblcirc_geom}, the union of two overlapping disks, with centers $P_{1}=\left(- \gamma,0\right)$ and $P_{2}=\left(\gamma,0\right)$, respectively, and each with radius $R$:
\begin{equation} \nonumber
\Omega =  \left\{\left(x,y\right) : |\left(x,y\right)-P_{1}| < R \right\}\cup \left\{\left(x,y\right) : |\left(x,y\right)-P_{2}|<R \right\}
\end{equation}
\noindent where $|\left(x,y\right)| = \sqrt{x^{2}+y^{2}}$ is the usual Euclidean vector norm, and $\gamma < R$.

This geometry is of interest due to, for example, its similarity to that of the radio frequency (RF) cavities used in the design of linear particle accelerators, and presents numerical difficulties due to the curvature of, and presence of corners in, the boundary. Our method avoids the staircase approximation used in typical finite difference methods to handle curved boundaries, which reduces accuracy to first order and may introduce spurious numerical diffraction.

\begin{figure}[ht!]
	\centering
	\subfigure[$t=0          $]{\includegraphics[trim={1cm 1cm 1cm 1cm},clip, width=0.3\textwidth]{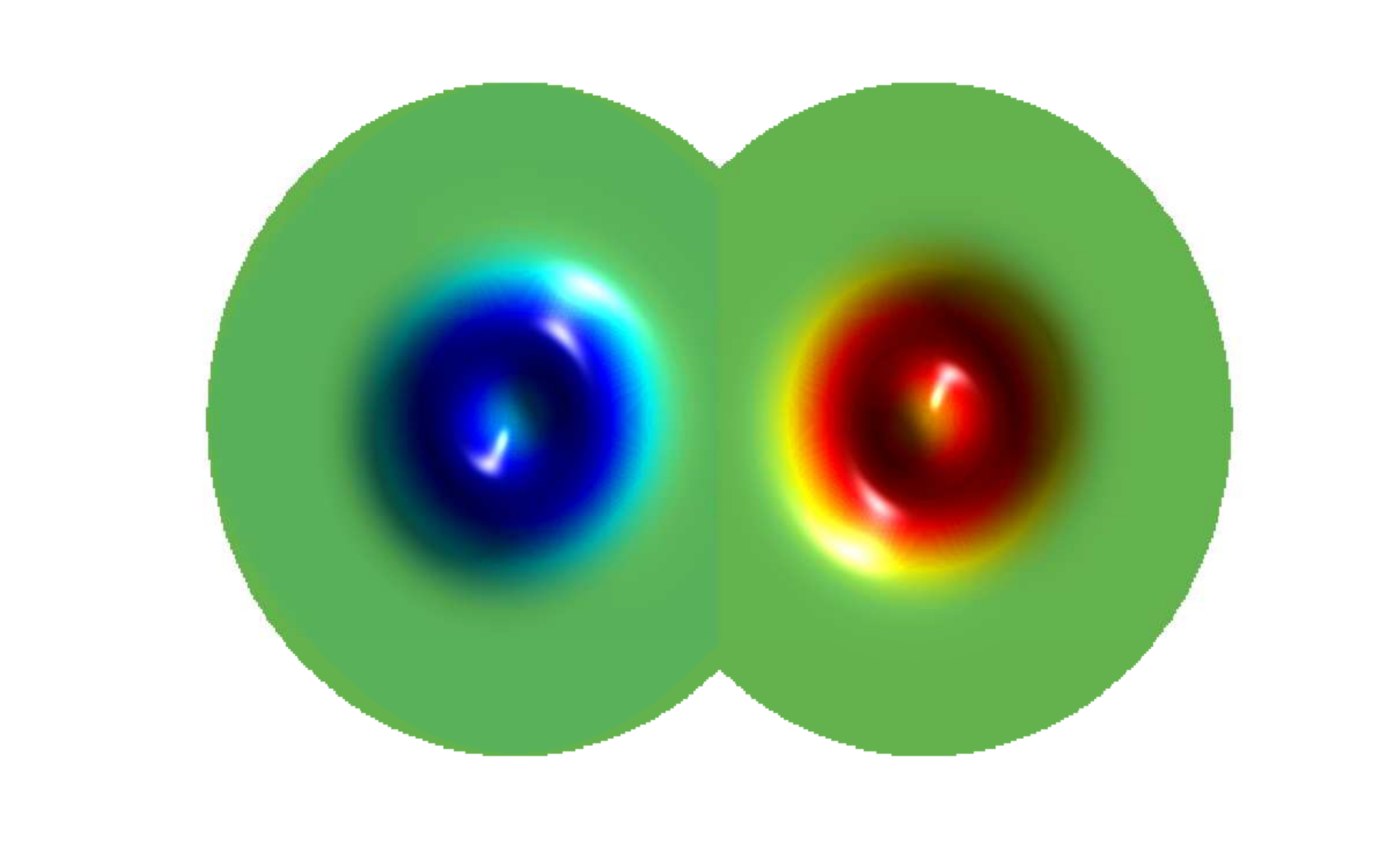}}
	\subfigure[$t=0.5145$]{\includegraphics[trim={1cm 1cm 1cm 1cm},clip, width=0.3\textwidth]{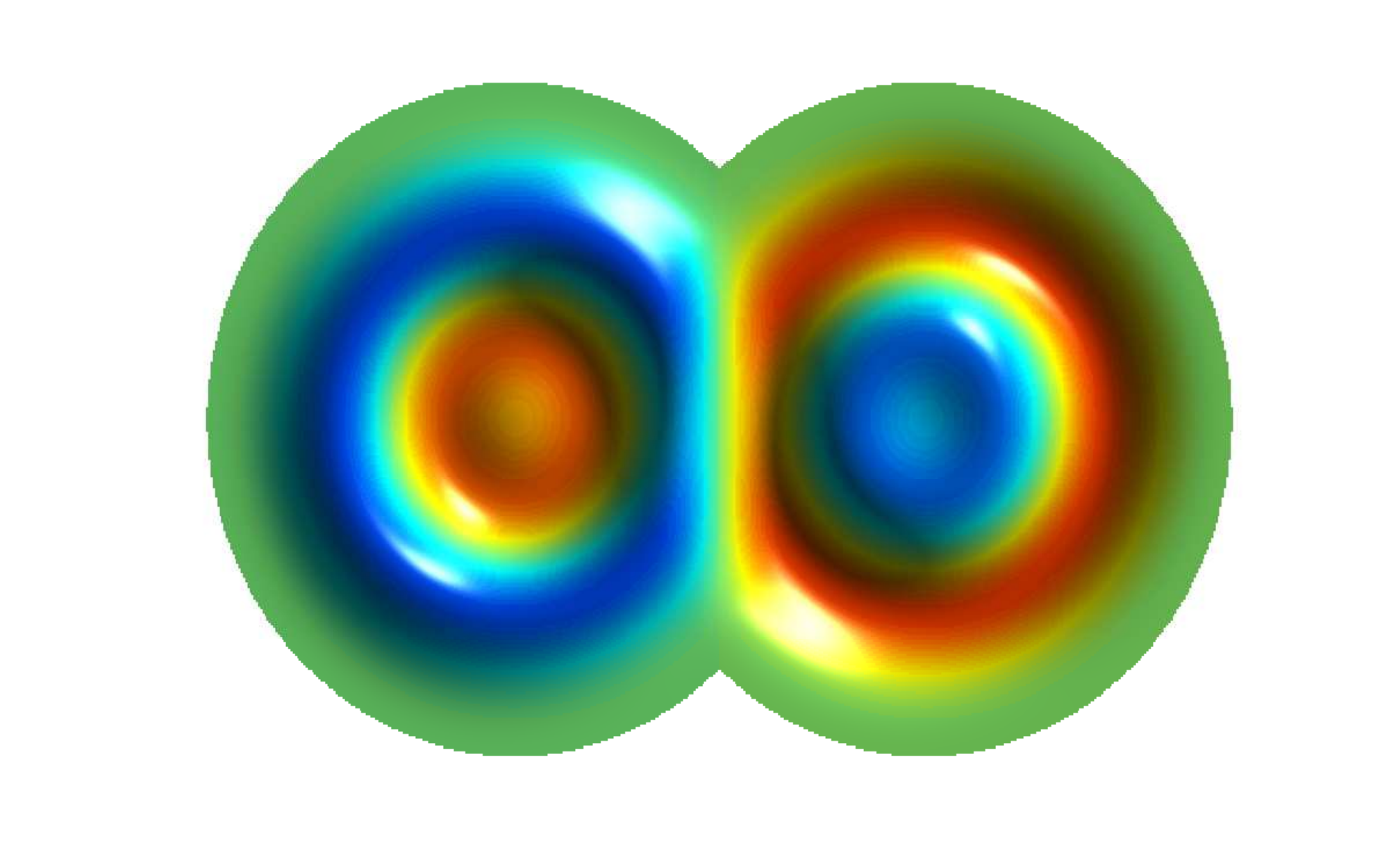}}
	\subfigure[$t=0.5145$]{\includegraphics[trim={1cm 1cm 1cm 1cm},clip, width=0.3\textwidth]{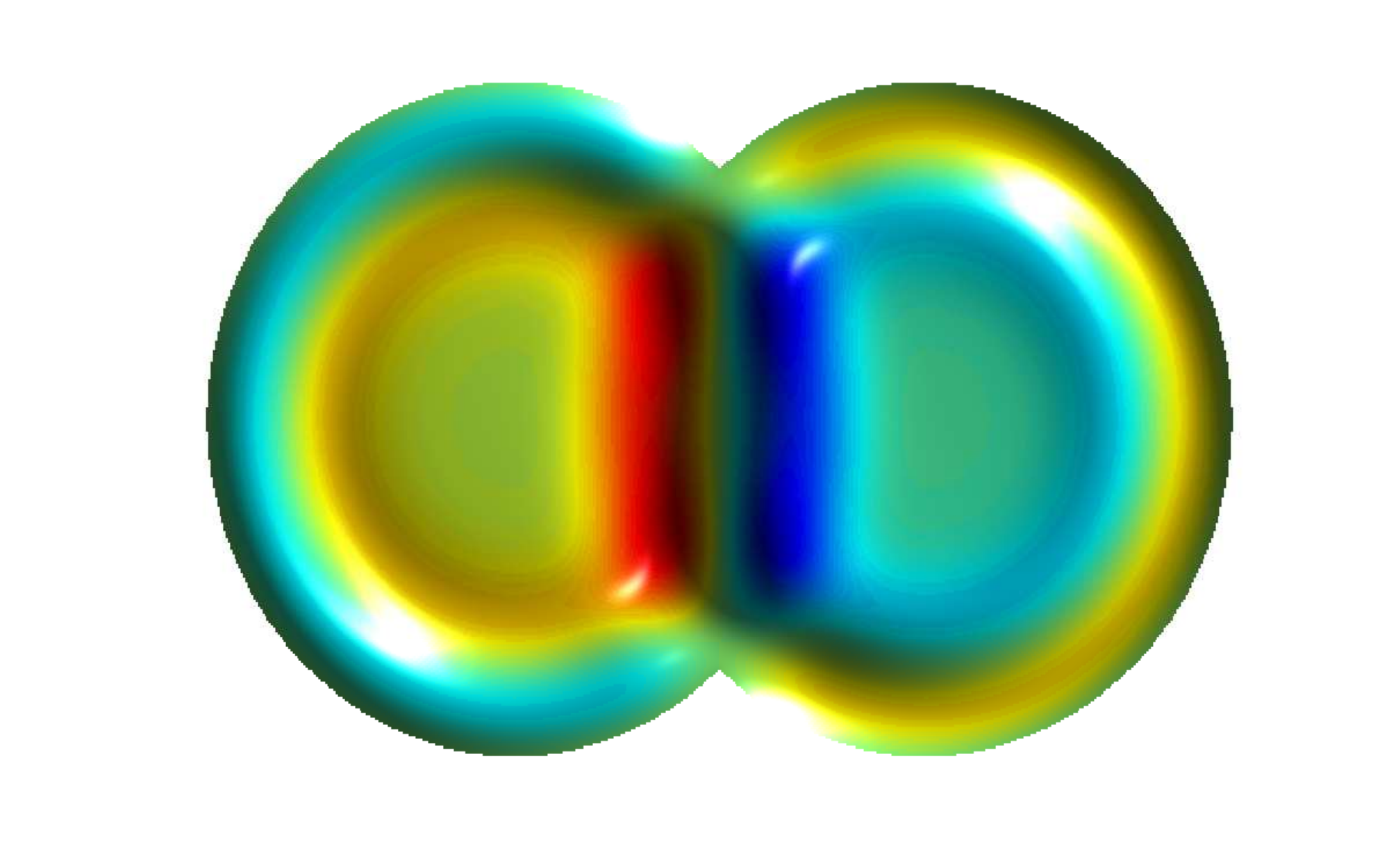}}
	\subfigure[$t=0.5145$]{\includegraphics[trim={1cm 1cm 1cm 1cm},clip, width=0.3\textwidth]{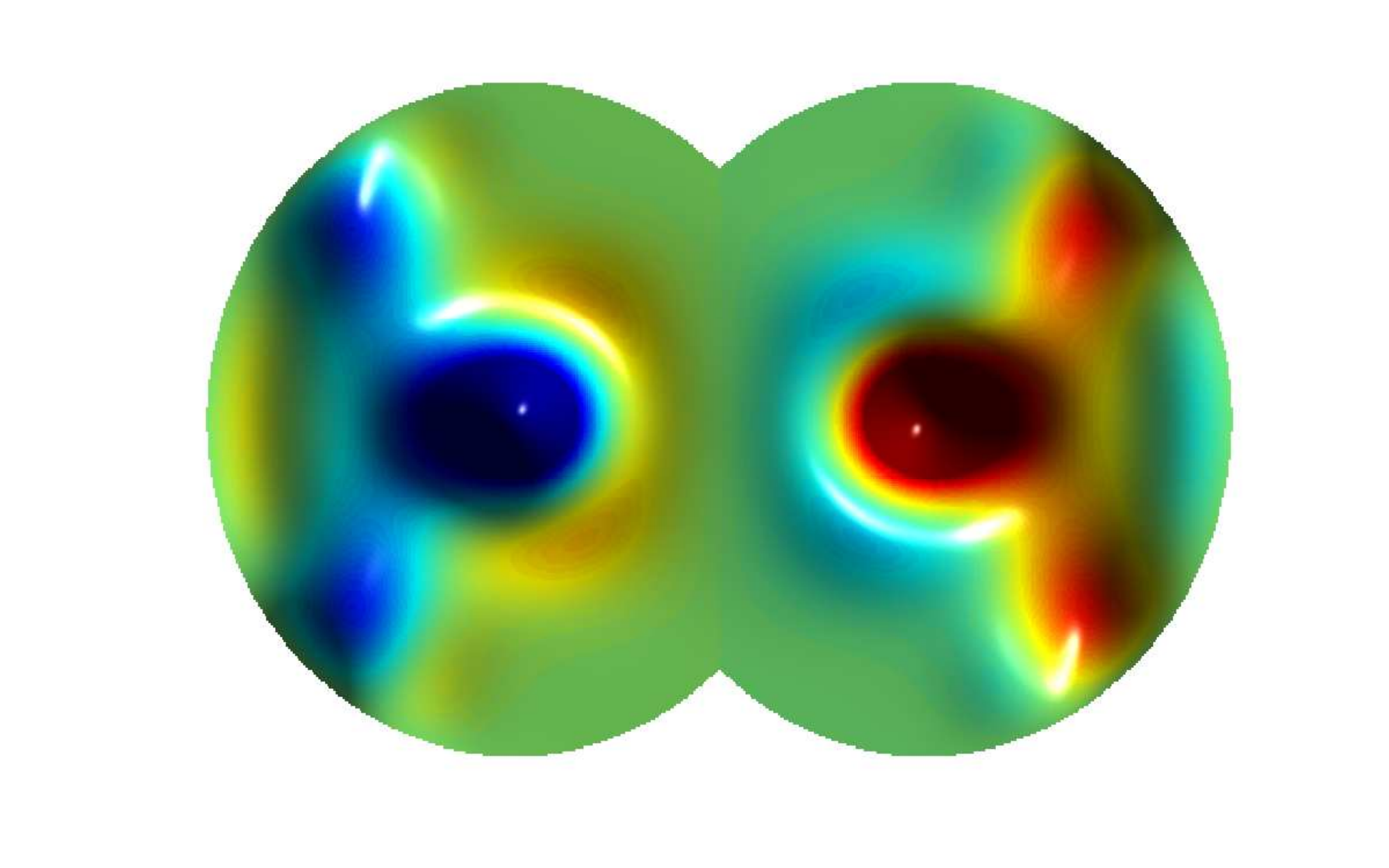}}
	\subfigure[$t=0.9135$]{\includegraphics[trim={1cm 1cm 1cm 1cm},clip, width=0.3\textwidth]{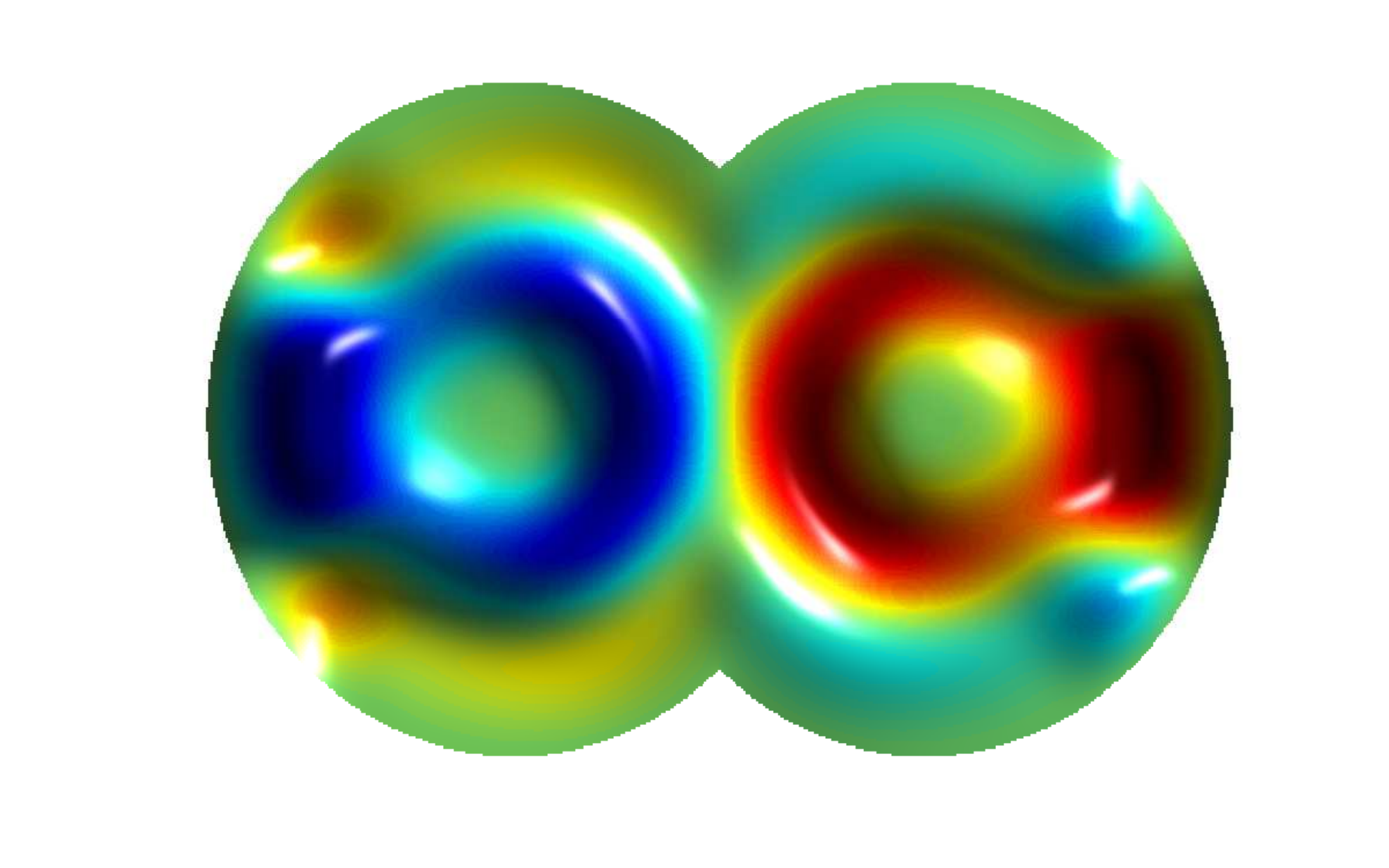}}
	\subfigure[$t=1           $]{\includegraphics[trim={1cm 1cm 1cm 1cm},clip, width=0.3\textwidth]{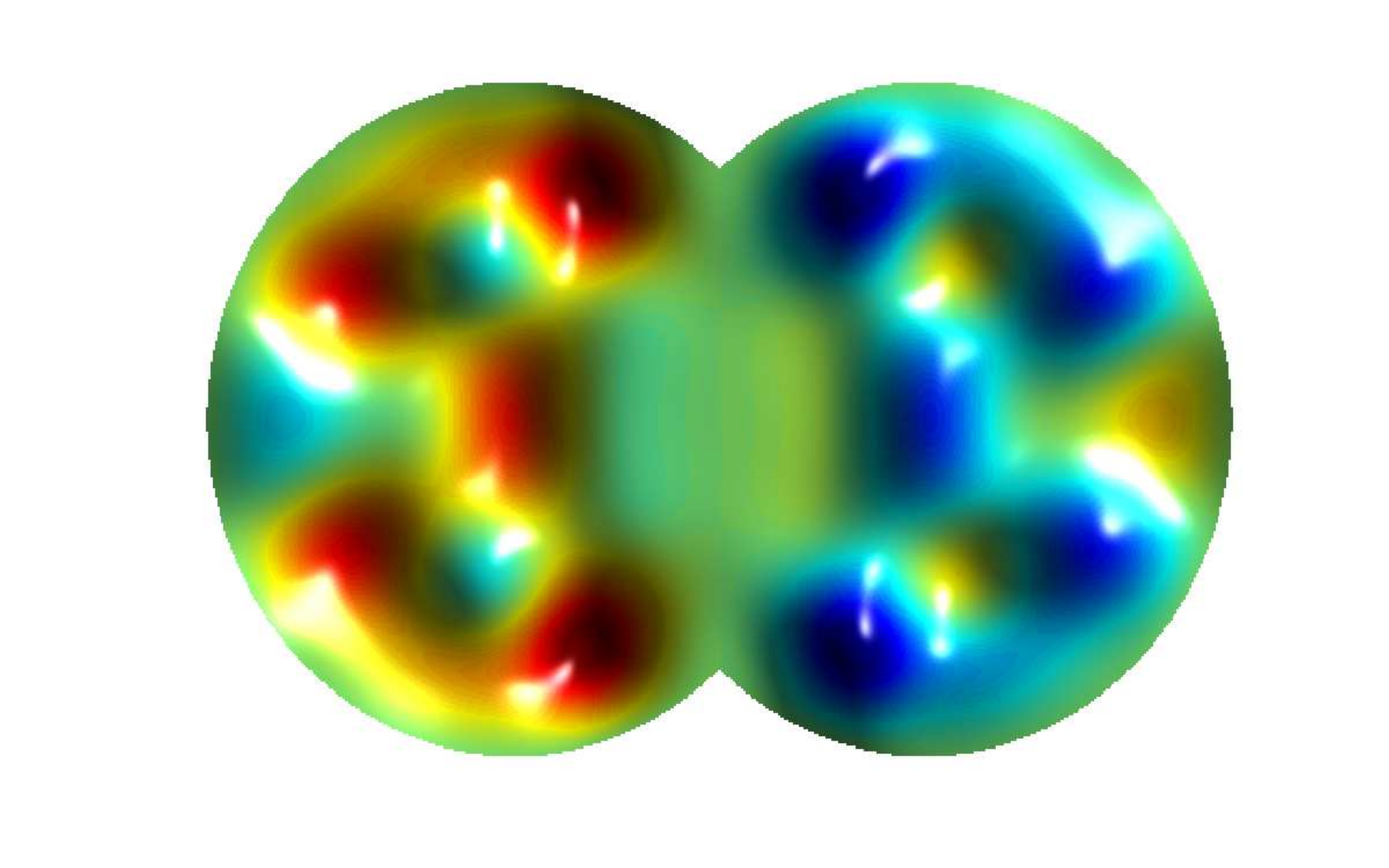}}
	\caption{Evolution of the double circle cavity problem.}
	\label{fig:dblcirc_plots}
\end{figure}

As initial conditions, we choose
\[
u\left(x,y,0\right)  = 
\begin{cases}
-\cos^{6}\left(\frac{\pi}{2}\left(\frac{|(x,y)-P_{1}|}{0.8 \gamma}\right)^{2}\right)	& |\left(x,y\right)-P_{1}| < 0.8 \gamma \\
\cos^{6}\left(\frac{\pi}{2}\left(\frac{|(x,y)-P_{2}|}{0.8 \gamma}\right)^{2}\right)	& |\left(x,y\right)-P_{2}| < 0.8 \gamma \\
0 & \text{otherwise}
\end{cases}
\]
and
\[
u_{t}\left(x,y,0\right)=0
\]
for $(x,y) \in \Omega$. Selected snapshots of the evolution are given in Figure \ref{fig:dblcirc_plots}, and the results of a refinement study are given in Table \ref{tab:refinement_dblcirc}. The discrete $L^{2}$ error was computed against a well-refined numerical reference solution ($\Delta x = \num{0.00021875}$); the error displayed in the table is the maximum over time steps with $t \in [0.28,0.29]$. For this example, $R=0.3$, $\gamma=0.2$, $c=1$, and the CFL is 2.

\begin{table}[htbp]
	\begin{centering}
		\begin{tabular}{|c|c|c||c|c|} \hline 
			$\Delta x$  		& $\Delta y$ 			& $\Delta t$			& $L^{2}$ error			& $L^{2}$ order	\\ \hline
			$\num{0.0070}$		& $\num{0.0043333333}$	& $\num{0.008666667}$	& $\num{0.006143688}$	& $-$			\\ \hline
			$\num{0.0035}$	& $\num{0.0021666667}$	& $\num{0.004333333}$	& $\num{0.001682923}$	& $1.8681$		\\ \hline
			$\num{0.00175}$	& $\num{0.0010833333}$	& $\num{0.00216666}$	& $\num{0.000435945}$	& $1.9488$		\\ \hline
			$\num{0.000875}$	& $\num{0.0005416667}$	& $\num{0.00108333}$	& $\num{0.000105150}$	& $2.0517$		\\ \hline
		\end{tabular}
		\caption{Refinement study for the double circle cavity with Dirichlet BC. For the numerical reference solution, $\Delta x = \num{.00021875}$, $\Delta y = \num{.00013542}$, and $\Delta t = \num{.00027083}$.}
		\label{tab:refinement_dblcirc}
	\end{centering}
\end{table}

\subsection{Symmetry on a Quarter Circle}

With the goal of testing the capabilities of our boundary conditions as well as circular geometry, we construct standing modes on a circular wave guide of radius $R$, in two different ways. First, we solve the Dirichlet problem, with initial conditions
\[
u(x,y,0) = J_0\left(z_{20}\frac{r}{R} \right), \quad u_t(x,y,0) = 0,
\]
and exact solution $u = J_0\left(z_{20}\frac{r}{R} \right)\cos\left(z_{20}\frac{ct}{R}\right)$, where $J_0$ is the Bessel function of order 0, and $z_{20} = 5.5218$ is the $2$-nd zero. Secondly, we use the symmetry of the mode to construct the solution restricted to the second quadrant, with homogeneous Neumann boundary conditions taken along the $x$ and $y$ axes.

\begin{figure}[hb]
	\centering
	\subfigure[$t = 0.25$]{\includegraphics[width = 0.24\textwidth]{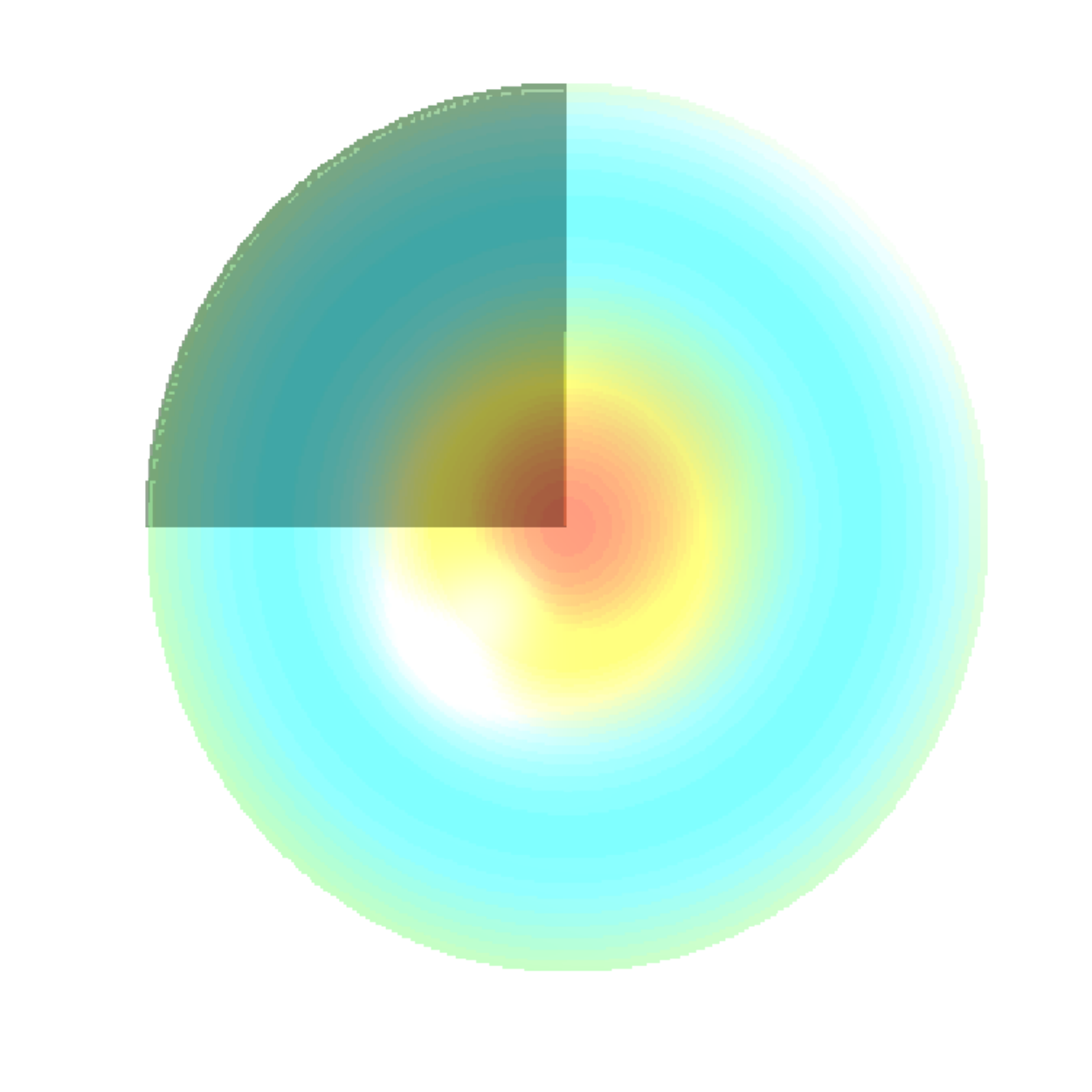}}
	\subfigure[$t = 0.50$]{\includegraphics[width = 0.24\textwidth]{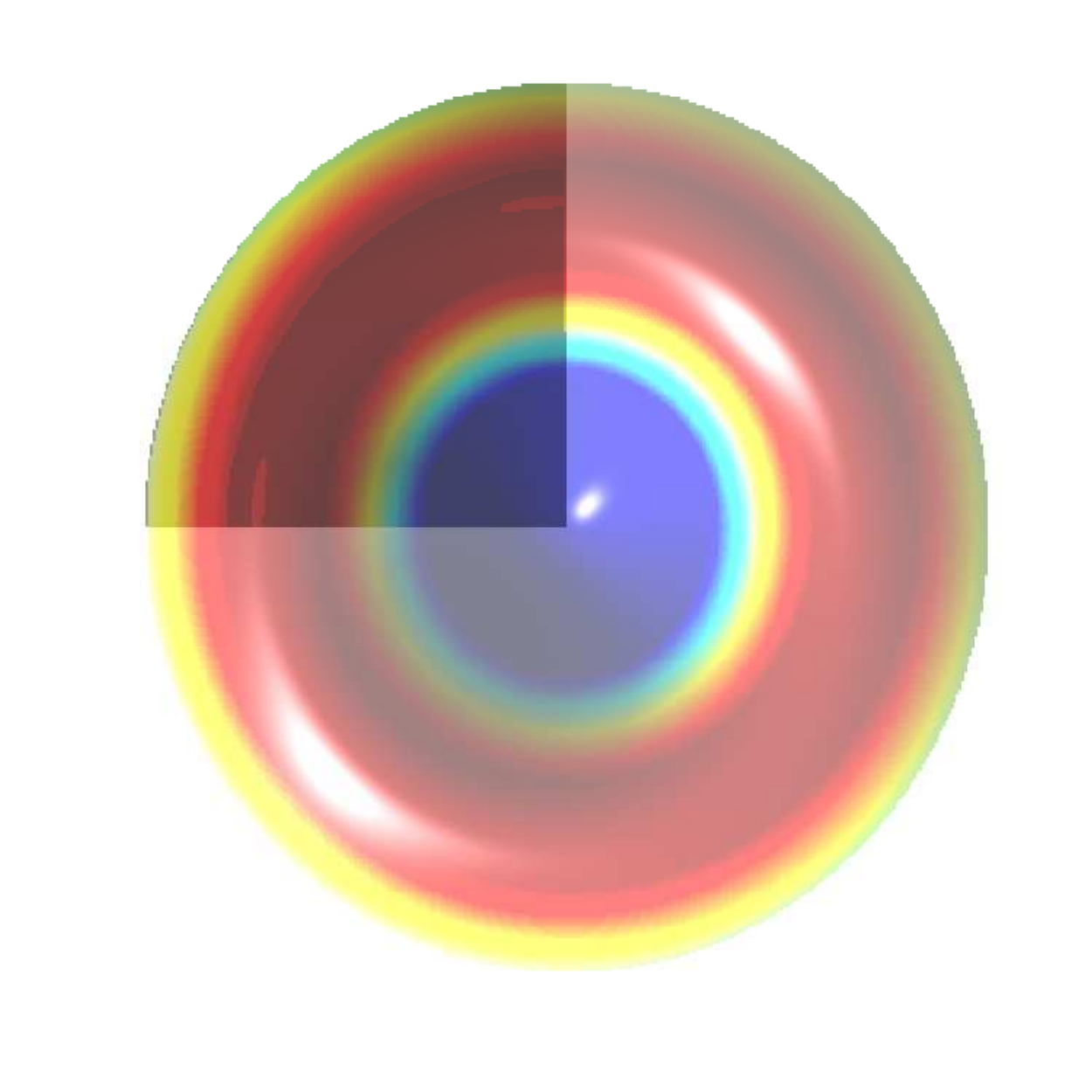}}
	\subfigure[$t = 0.75$]{\includegraphics[width = 0.24\textwidth]{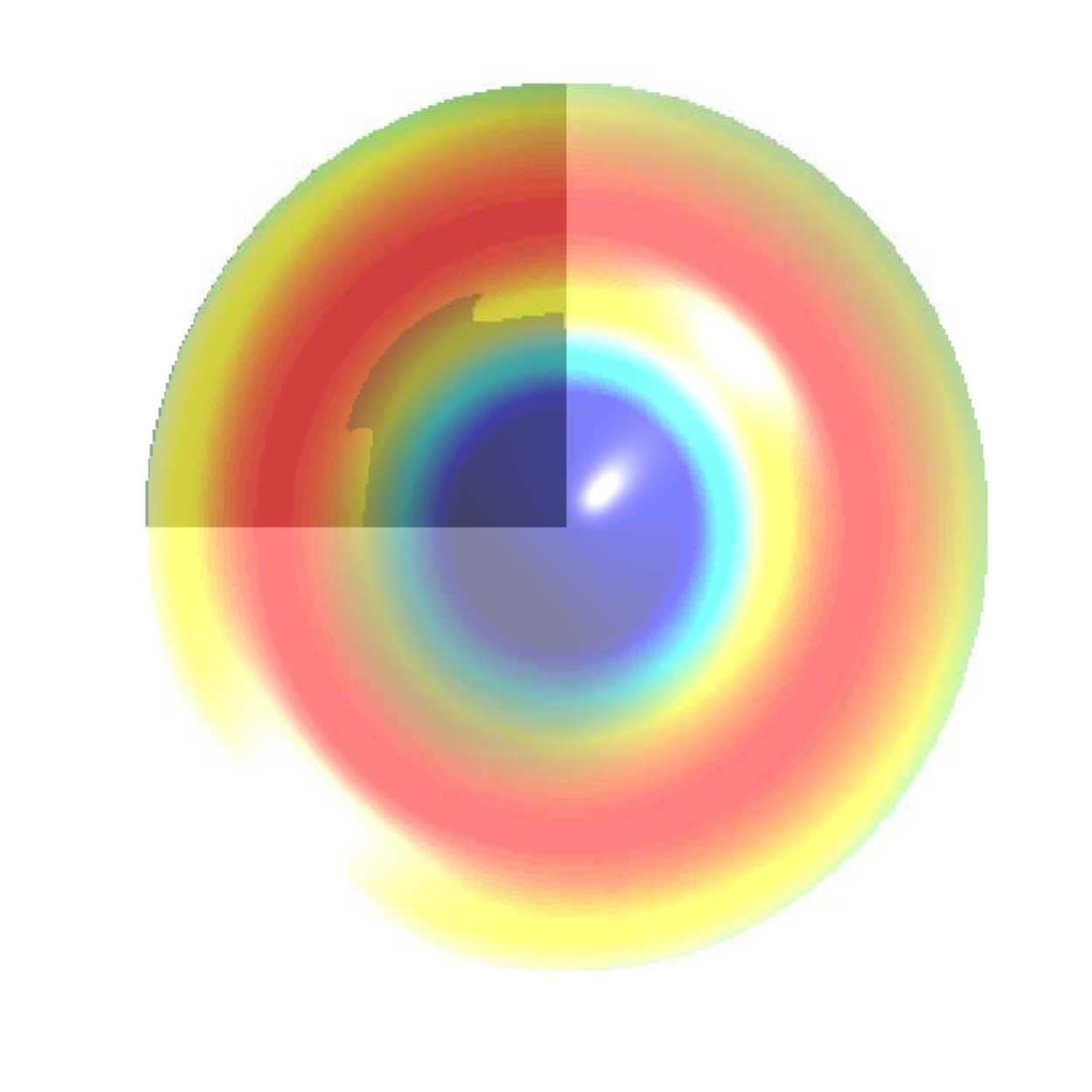}}
	\subfigure[$t = 1.00$]{\includegraphics[width = 0.24\textwidth]{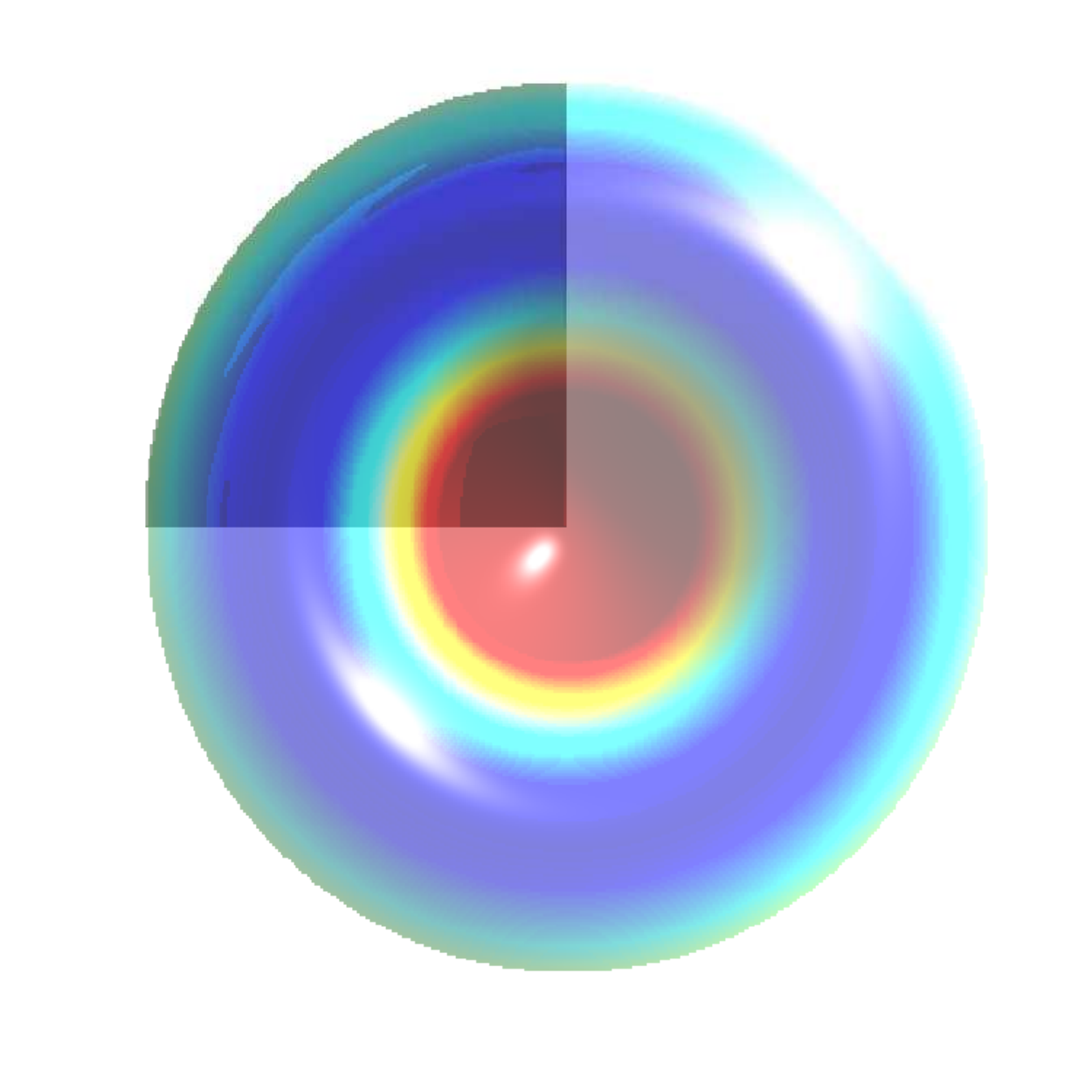}}
	\caption{Two separate numerical constructions of a Bessel mode are superimposed, demonstrating that the solution on the quarter circle using Neumann boundary conditions is equivalent to that of the full circle.}
	\label{fig:Bessel_Quarter}
\end{figure}

In both cases, the solution converges to second order. An overlay of the two are shown in Figure \ref{fig:Bessel_Quarter}, demonstrating the close agreement.

\subsection{Bessel Mode with Neumann Boundary Conditions}

\begin{figure} 
	\begin{center}
		\includegraphics[width = \textwidth]{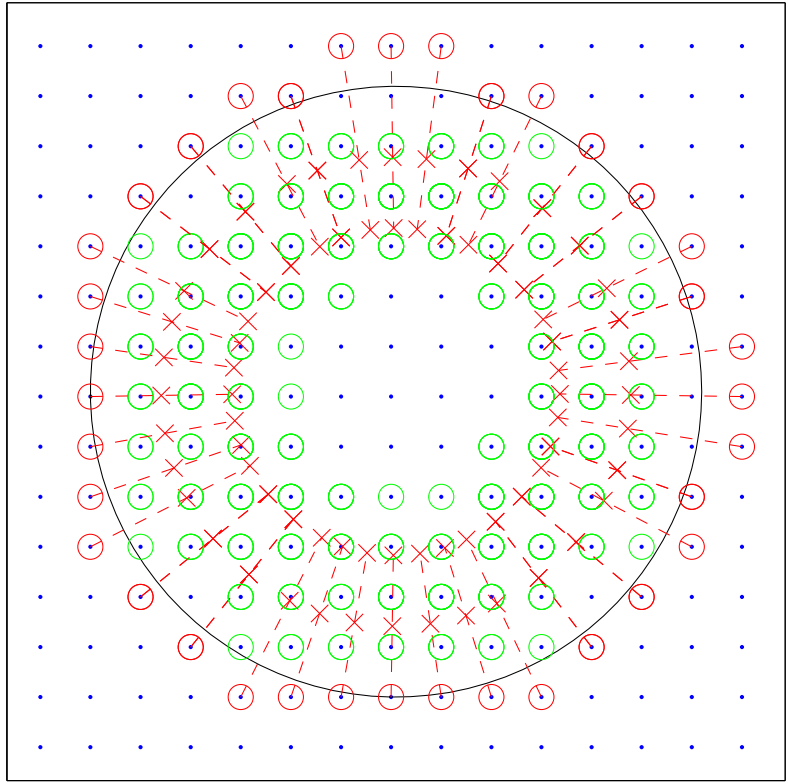}
	\end{center}
	\caption{An example of the embedded boundary grid used. The red circled exterior grid points are the endpoints where a value is to be calculated via the interpolation procedure. The red crosses are the points where values are imposed on quadratic boundary interpolant along the normal direction (red dashed line). Values for the bilinear interpolants are supplied from the green circled interior grid points.}
	\label{fig:eb_grid_plot}
\end{figure}

\begin{figure}
	\begin{center}
		\includegraphics[width = \textwidth]{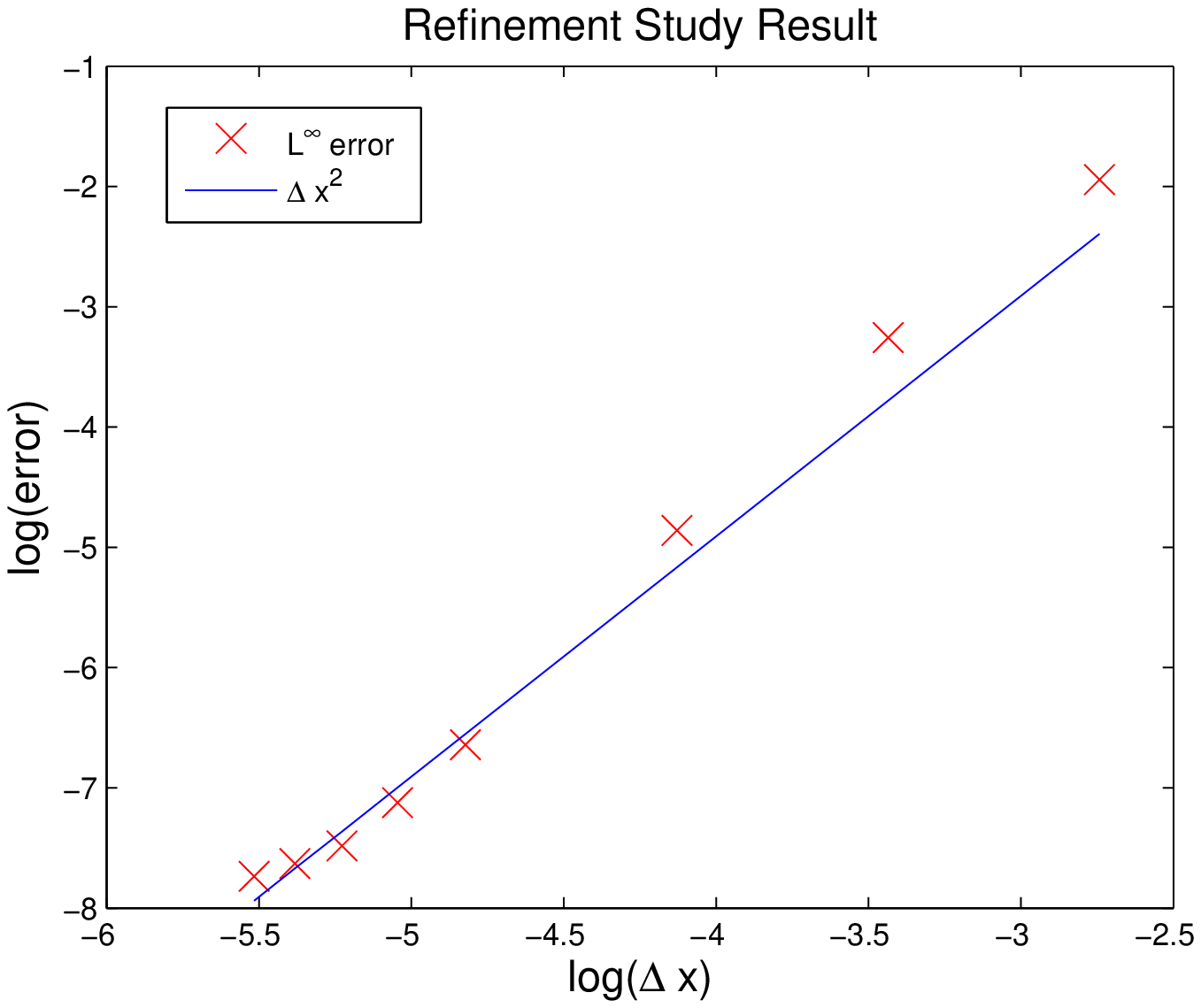}
	\end{center}
	\caption{Refinement study for the Bessel mode in a circular domain with fixed CFL number of 2. Using quadratic boundary interpolant with bilinear interior interpolant.}
	\label{fig:eb_refinement_study}
\end{figure}

Here we present a numerical example of the embedded boundary method for homogeneous Neumann boundary conditions given in Section \ref{sec:embedded_neumann_bc}. We apply the method to a circular domain, for which analytical solutions exist. We consider a radially-symmetric Bessel mode with homogeneous Neumann boundary conditions, with an analytic solution given by

\begin{align}
	u(r,t) &= J_{0}\left(Z_{0}\frac{r}{R}\right)\cos\left(Z_{0}\frac{ct}{R}\right),
\end{align}

where $J_{0}$ is the Bessel function of the first kind of order 0, $r = \sqrt{x^{2}+y^{2}}$, $R$ is the radius of the domain, and $Z_{0}\approx3.8317$ is the smallest nonzero root of $J_{0}'$ (so that $\frac{\partial u}{\partial \mathbf{n}}(R,t) = \frac{\partial u}{\partial r}(R,t) = \frac{Z_{0}}{R}J_{0}' \left(Z_{0}\right)\cos\left(Z_{0}\frac{ct}{R}\right) =0$). In this example, we take radius $R = \pi/2$ and wave speed $c=1$. An example of the embedded boundary grid used is given in Figure \ref{fig:eb_grid_plot}.  We perform a refinement study with a fixed CFL number of 2, with the results in Figure \ref{fig:eb_refinement_study} indicating the expected second-order convergence. We set the iteration tolerance to $10^{-15}$, and we see convergence of the boundary correction iteration in fewer than 40 iterations. We note some oscillation of the $L^{\infty}$ error about the line giving second-order accuracy, which we believe to be due to the grid points moving with respect to the boundary through the refinement, causing some variation in the maximum error.

\subsection{Periodic Slit Diffraction Grating}

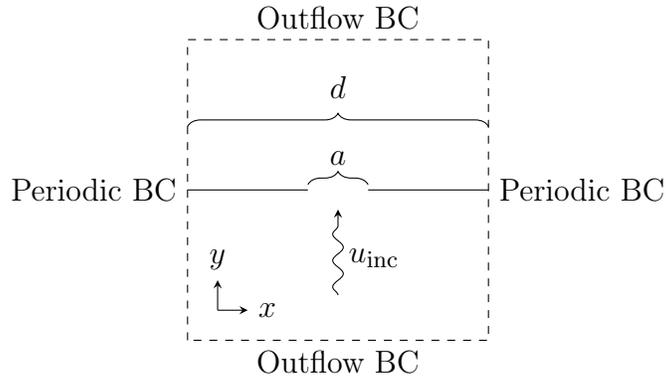
\begin{figure}[H]
	\centering
	\begin{tikzpicture}[scale=2,>=stealth]
	
	\draw (-1,0)--(-0.2,0);
	\draw (0.2,0)--(1,0);
	
	\draw[->] (-0.8,-0.8) -- (-0.6,-0.8)
	node[right] {$x$};
	\draw[->] (-0.8,-0.8) -- (-0.8,-0.6)
	node[above]{$y$};
	
	\draw [decorate,decoration={brace,amplitude=5pt},xshift=0pt,yshift=1pt]
	(-0.2,0) -- (0.2,0) node [black,midway,xshift=0,yshift=10pt] {$a$};
	\draw [decorate,decoration={brace,amplitude=5pt},xshift=0pt,yshift=12pt]
	(-1,0) -- (1,0) node [black,midway,xshift=0,yshift=15pt] {$d$};
	
	\draw[->,decorate,decoration={snake,amplitude=2,post length=.08}] 		
	(0,-.7) -- (0,-.2)
	node [right,align=center,midway]
	{
		{$u_{\mbox{\scriptsize inc}}$}
	};
	
	\draw[dashed] (-1,-1)--(-1,1)--(1,1)--(1,-1)--cycle;
	
	\draw[anchor=north] (0,-1) node {Outflow BC};
	\draw[anchor=south] (0,1) node {Outflow BC};
	\draw[anchor=west] (1,0) node {Periodic BC};
	\draw[anchor=east] (-1,0) node {Periodic BC};
	
	\end{tikzpicture}
	\caption{Periodic slit diffraction grating geometry}
	\label{fig:slit_geom}
\end{figure}

In this example, we apply our method to model an infinite, periodic diffraction grating under an incident plane wave.  Diffraction gratings are periodic structures used in optics to separate different wavelengths of light, much like a prism. The high resolution that can be achieved with diffraction gratings makes them useful in spectroscopy, for example, in the determination of atomic and molecular spectra. Our numerical experiment, depicted in Figure \ref{fig:slit_geom}, demonstrates the use of our method with multiple boundary conditions and nontrivial geometry in a single simulation to capture complex wave phenomena.

An idealized slit diffraction grating consists of a reflecting screen of vanishing thickness, with open slits of aperture width $a$, spaced distance $d$ apart, measured from the end of one slit to the beginning of the next (that is, the periodicity of the grating is $d$). We impose an incident plane wave of the form $u_{inc}(x,y,t)=\cos{\left(\omega t + ky\right)}$, where $k = 2 \pi / a$ and $\omega = k/c$, where $c$ is the wave speed. Periodic BCs at $x=\pm d/2$ (determining the periodicity of the grating), and homogeneous Dirichlet BCs are imposed at the screen. We also test outflow boundary conditions in multiple dimensions, which are imposed at $y=\pm L_{y}/2$.

\begin{figure}[htbp!]
	\centering
	\subfigure[$t=0.31$]{\includegraphics[trim={1cm 1cm 1cm 1cm},clip, width = .24\textwidth]{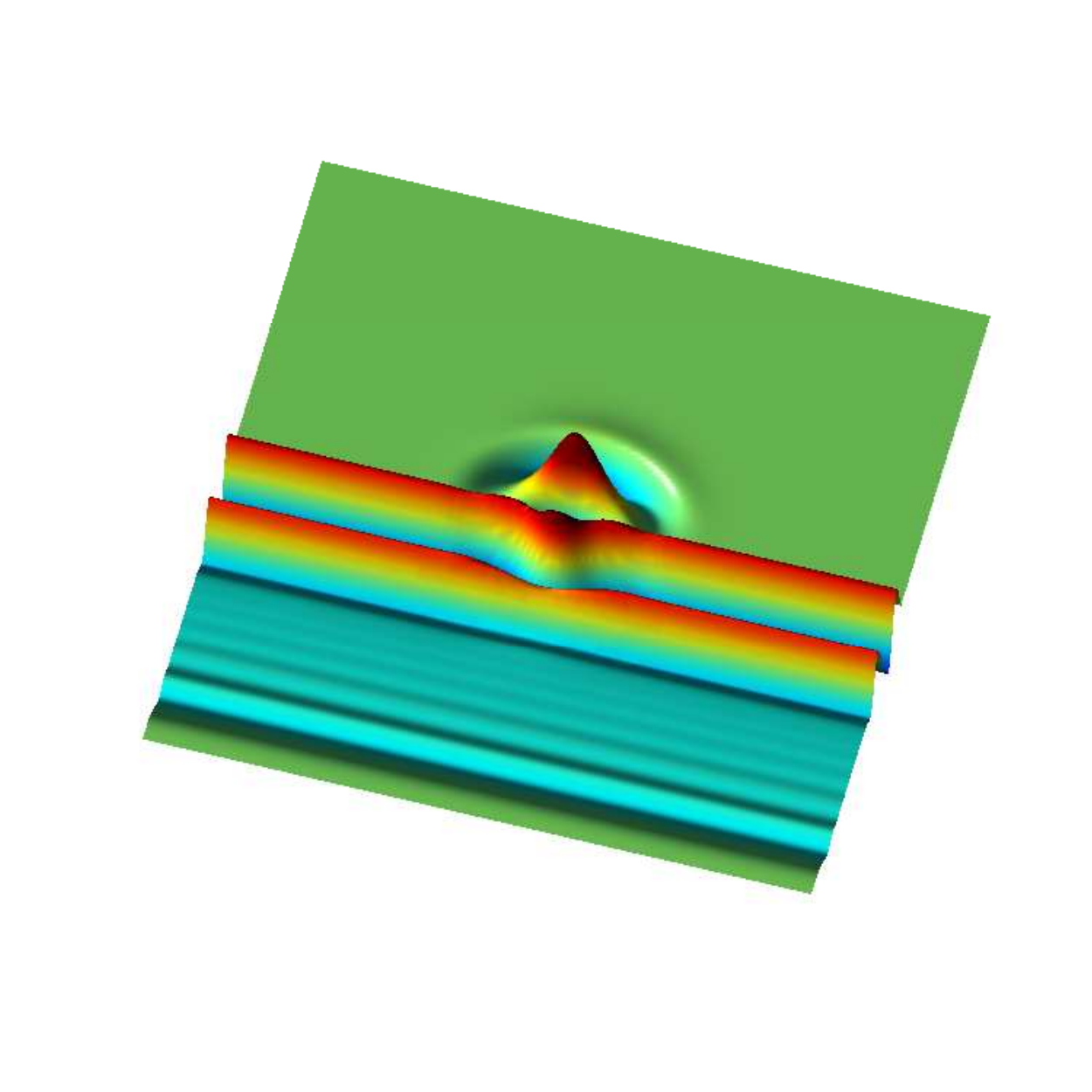}}
	\subfigure[$t=0.51$]{\includegraphics[trim={1cm 1cm 1cm 1cm},clip, width = .24\textwidth]{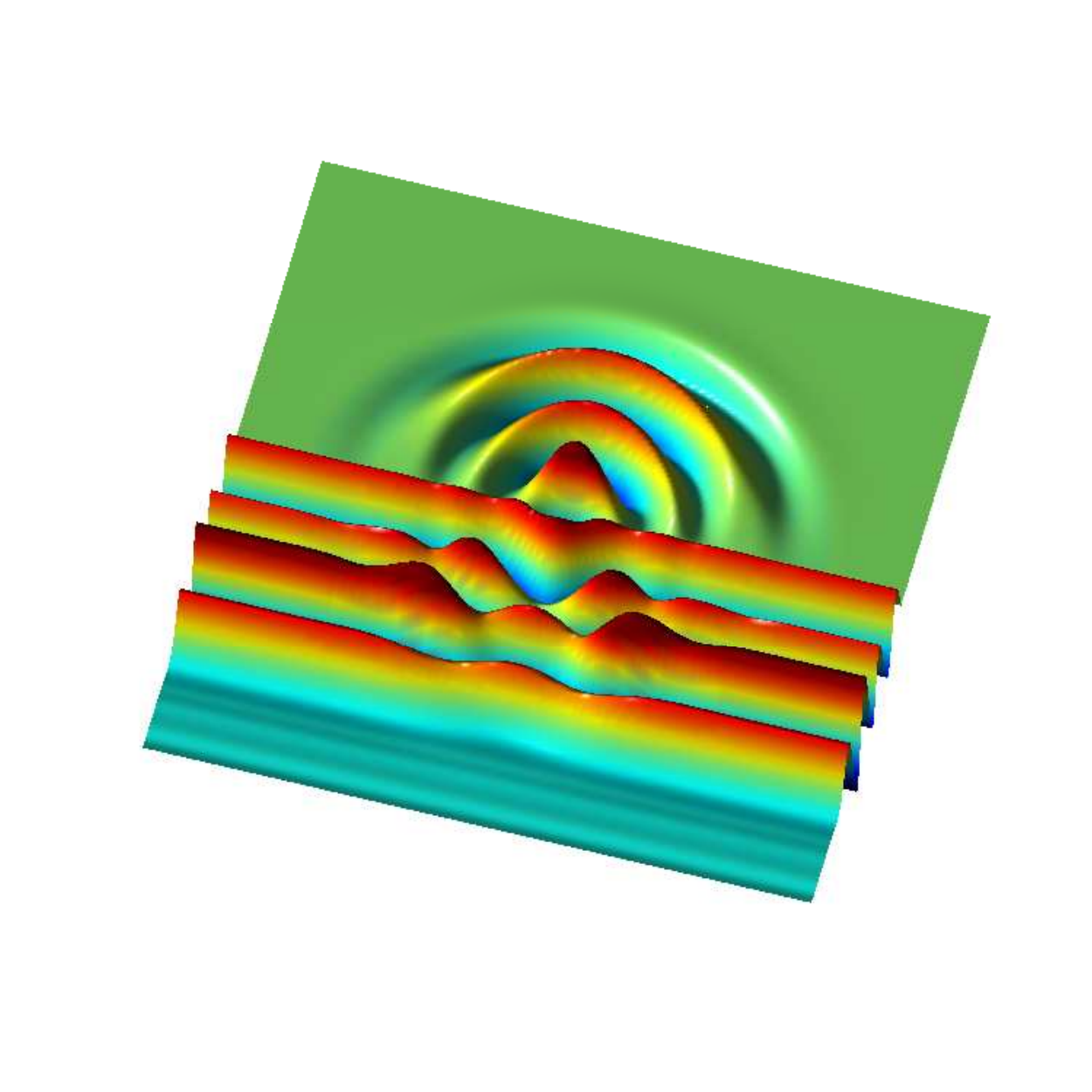}}
	\subfigure[$t=1.01$]{\includegraphics[trim={1cm 1cm 1cm 1cm},clip, width = .24\textwidth]{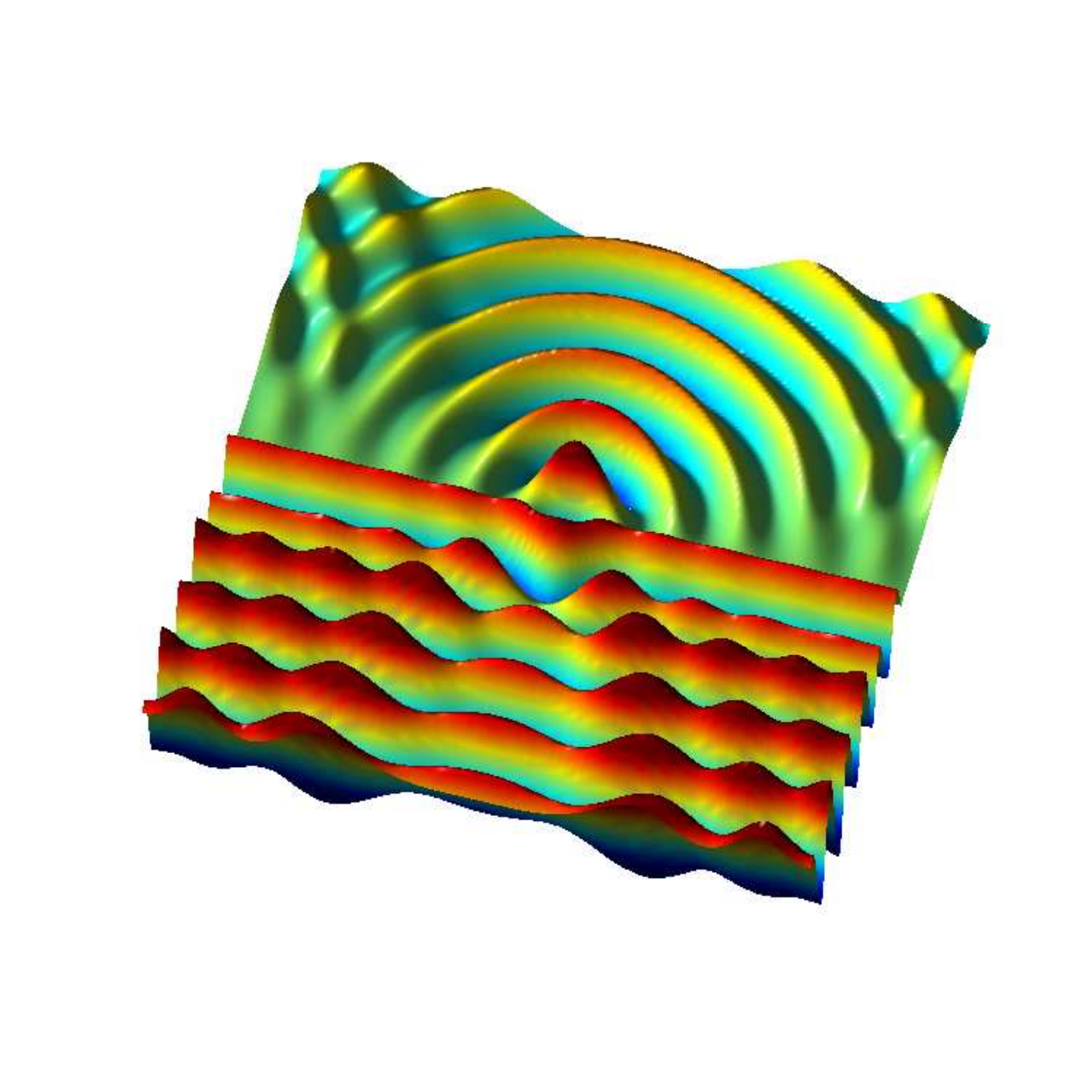}}
	\subfigure[$t=2.01$]{\includegraphics[trim={1cm 1cm 1cm 1cm},clip, width = .24\textwidth]{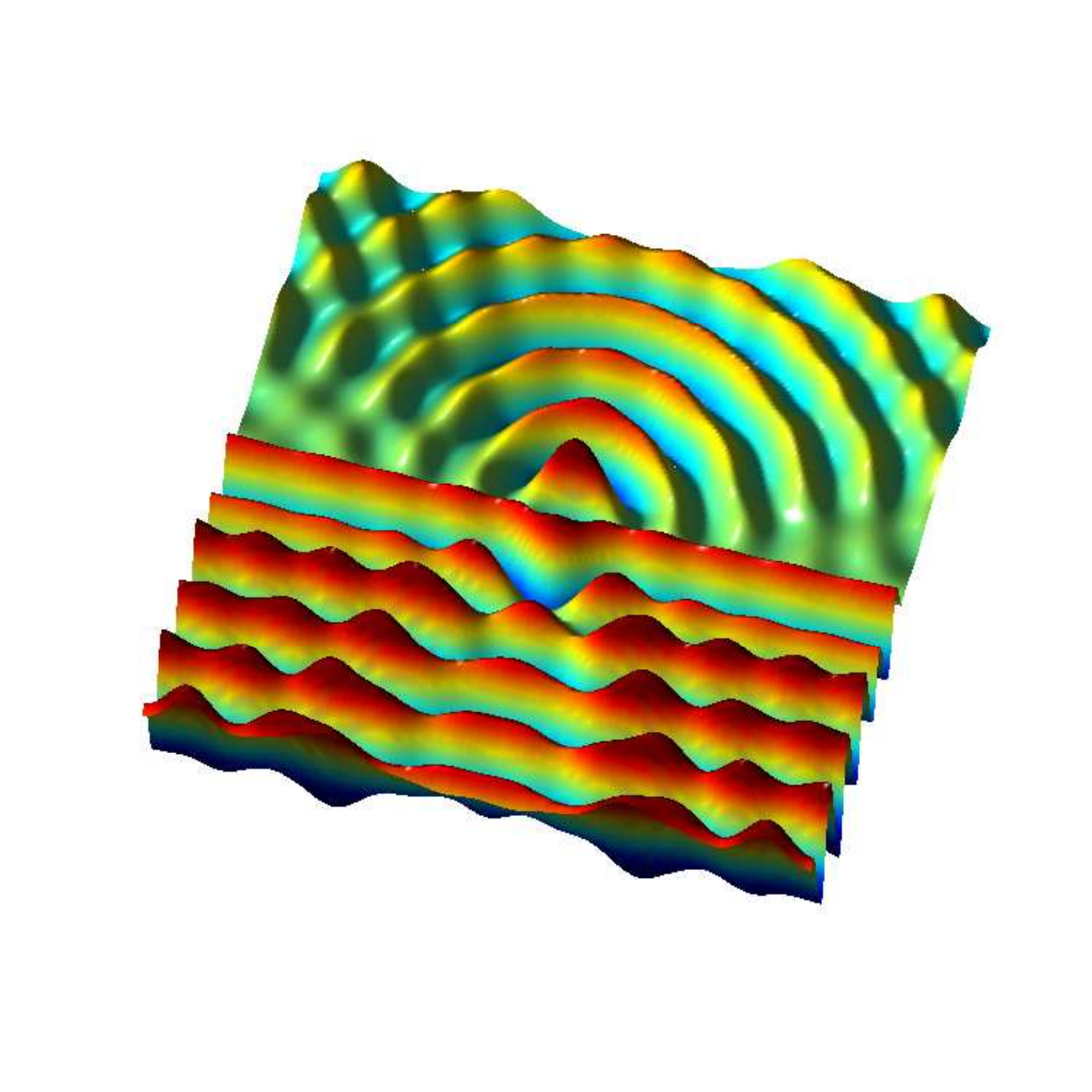}}
	\caption{Evolution of the slit diffraction grating problem, with aperture width $a=0.1$, grating periodicity $L_{y}=d=1$, and wave speed $c=1$. The CFL is fixed at 2.}
	\label{fig:slit_plots}
\end{figure}

In Figure \ref{fig:slit_plots}, we observe the time evolution of the incident plane wave passing through the aperture, and the resulting interference patterns as the diffracted wave propagates across the periodic boundaries. The outflow boundary conditions allow the waves to propagate outside the domain. While a rigorous analysis of the efficacy of our outflow BCs is the subject of future work, the results look quite reasonable, as no spurious reflections are seen at the artificial boundaries.

\subsection{Point Sources and Outflow Boundary Conditions}
\begin{figure}[htbp!]
	\centering
	\subfigure[$t=0.05$]{\includegraphics[width = .24\textwidth]{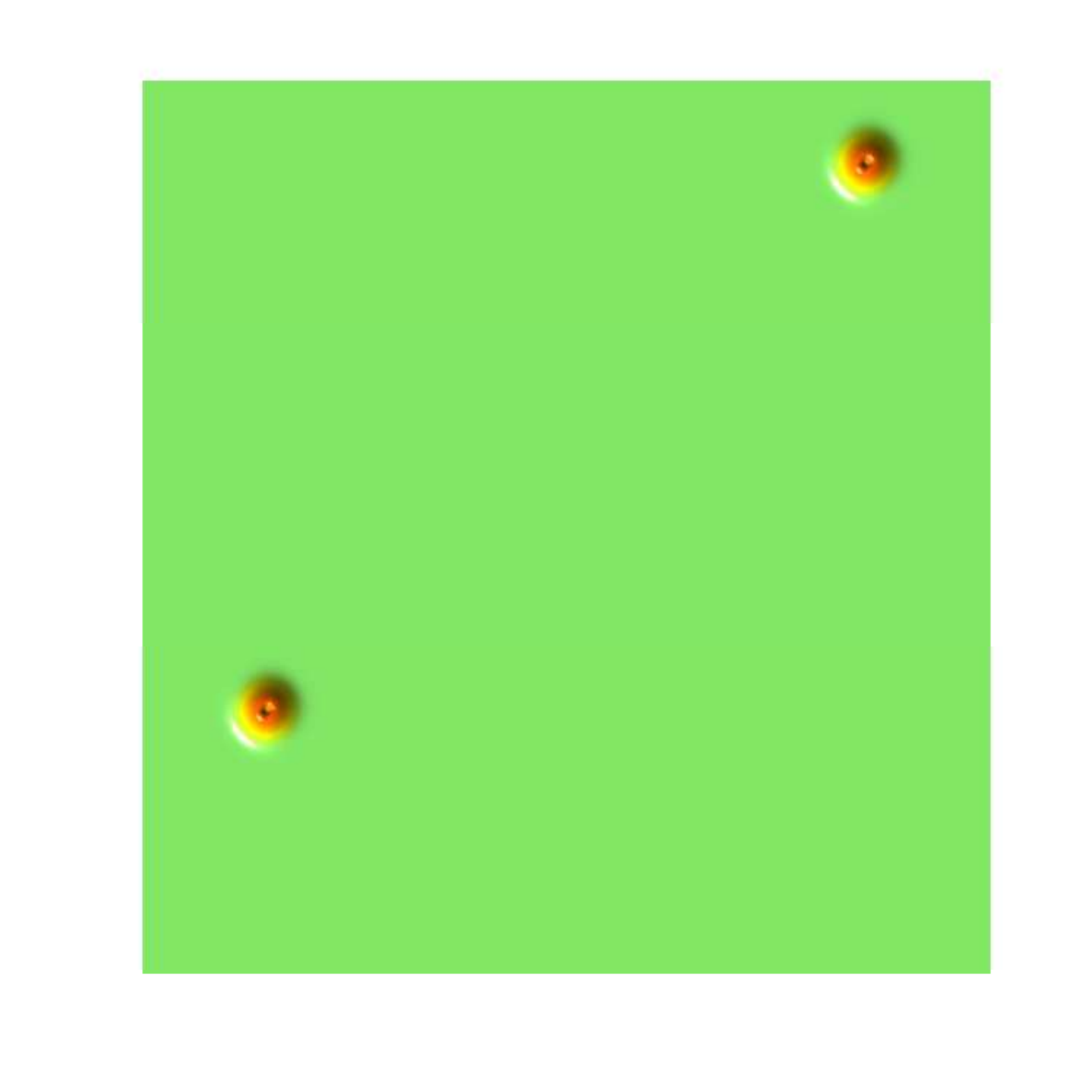}}
	\subfigure[$t=0.40$]{\includegraphics[width = .24\textwidth]{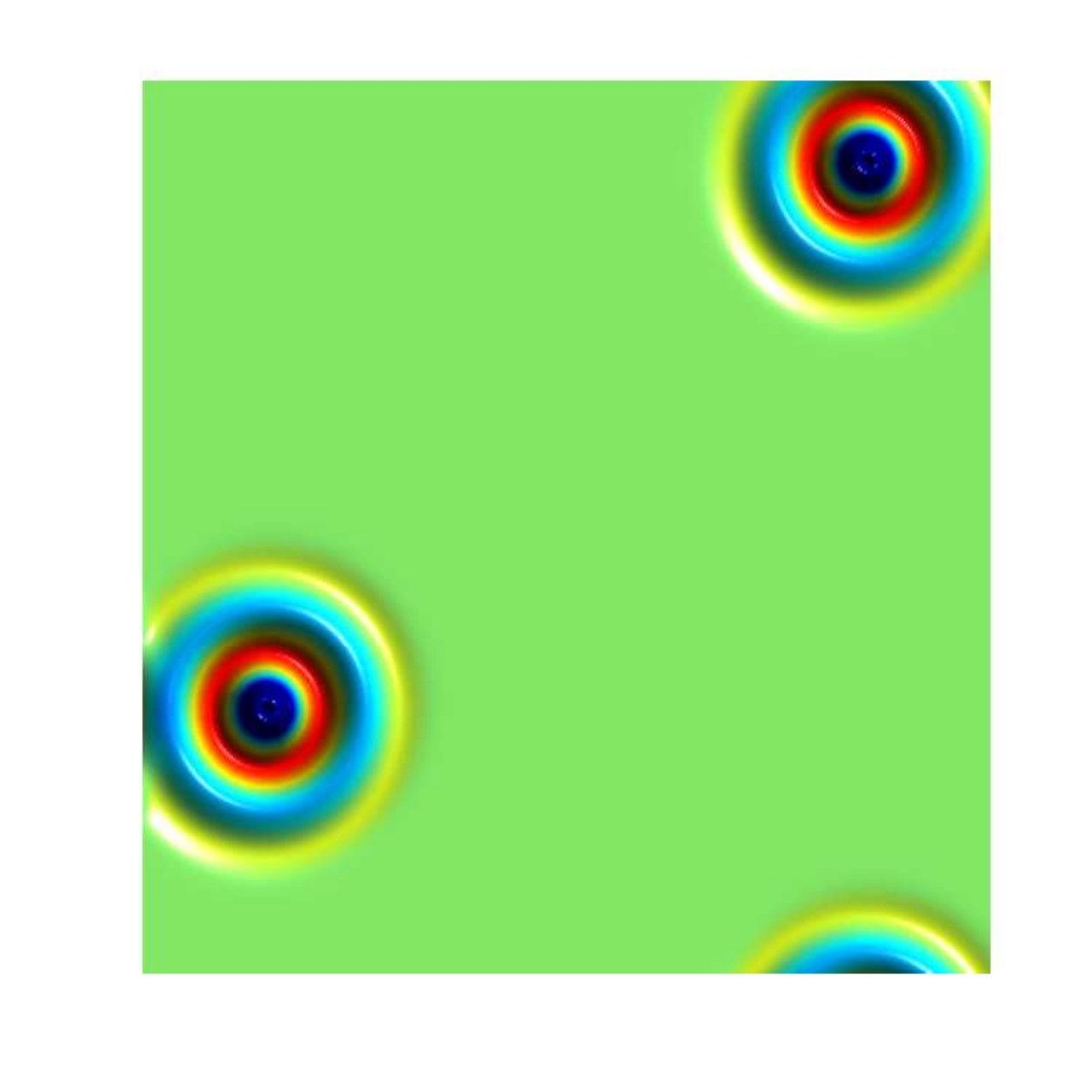}}
	\subfigure[$t=0.75$]{\includegraphics[width = .24\textwidth]{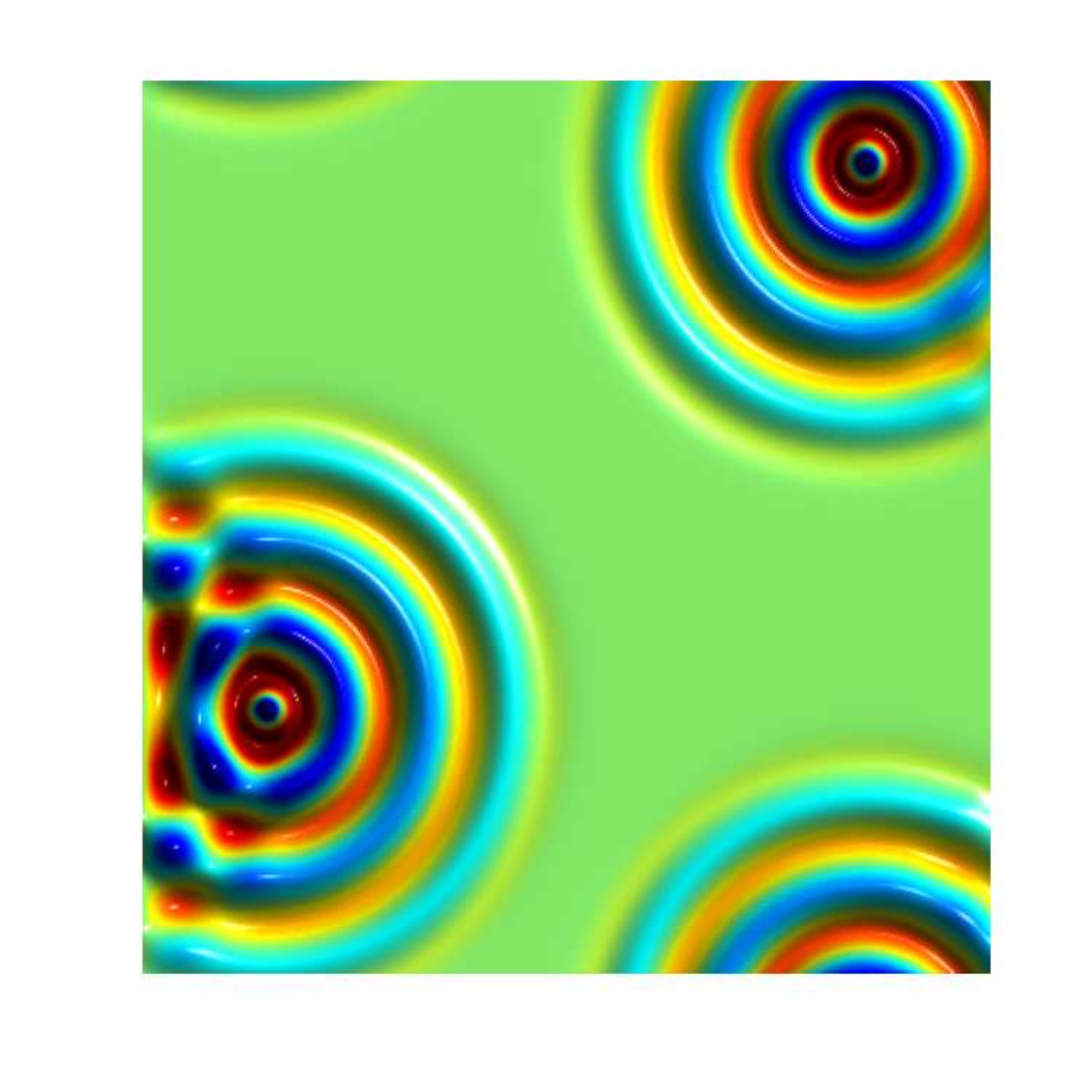}}
	\subfigure[$t=1.00$]{\includegraphics[width = .24\textwidth]{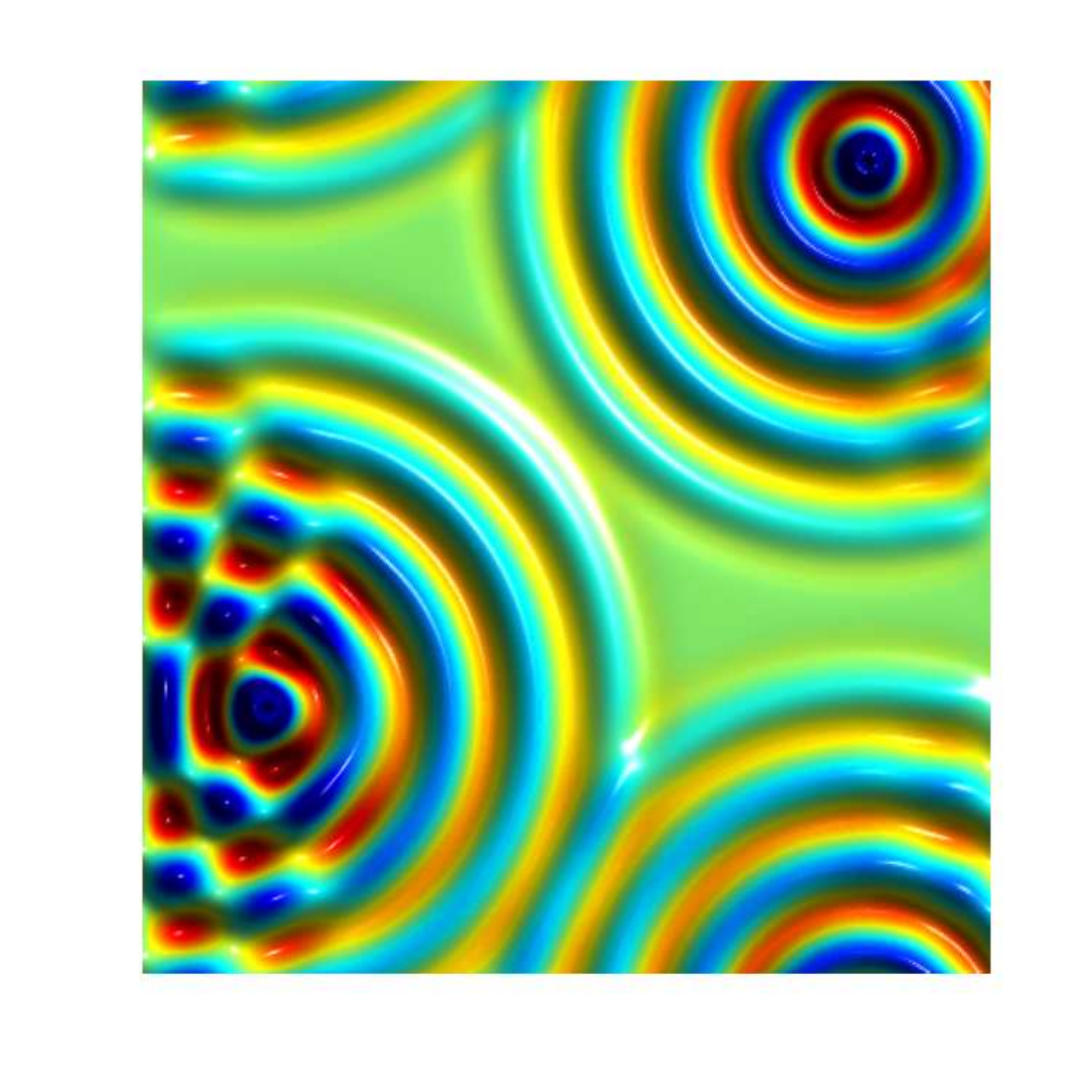}}
	\caption{Point Sources emanating through the domain. Outflow boundary conditions are prescribed at the right, Dirichlet boundary conditions on the left, and the top and bottom edges are periodic.}
	\label{fig:Source}
\end{figure}

As a final example, we illustrate some of the more interesting features of our solver. We launch two point sources (from points not on the mesh), and employ Dirichlet, periodic and outflow boundary conditions along the edges of the domain. As can be seen in Figure \ref{fig:Source}, the point sources propagate perfectly, despite not being placed at grid points.

\section{Conclusion}
\label{sec:conclusion}

In this paper we have presented a fast, A-stable, second order scheme for solving the wave equation. Using the method of lines transpose (MOL$^T$), the solution can be interpreted in the semi-discrete sense as a boundary integral solution, posed as a convolution against an exponential kernel. We have exploited this fact to develop a matrix-free, $O(N)$ fast spatial convolution algorithm, capable of embedding boundary points in a regular Cartesian mesh, without affecting the accuracy or stability. 

In addition to demonstrating second order accuracy for wave propagation in a variety of non-Cartesian geometries, with time steps in excess of the usual CFL restriction, we have also developed a novel method for implementing outflow boundary conditions, as well as methods for launching waves, using soft sources, from points located at arbitrary locations (e.g., not located at a mesh point) inside the domain, which is of interest in particle-wave simulations such as those required in studying plasma. Several topics warrant further investigation. We are developing a domain decomposition approach to multi-core computing with our implicit wave solver based on the properties of the exponential kernel, where the subdomains require only pointwise data communication from adjoining edges. Future work will investigate the implementation of higher-order accurate boundary conditions for complex boundary geometries and outflow boundary conditions, and further, apply these methods to Maxwell's equations both in electromagnetic scattering and plasma physics problems.

\appendix
\section{Summary Table For Second-Order Wave Solver}
\label{sec:wave_summary}

\begin{center}
	\renewcommand{\arraystretch}{1}
	\small \begin{tabular}{|c|c|}
		\hline
		Wave Equation & Dispersive Scheme, $\alpha = \frac{2}{c \Delta t}$: \\
		$\frac{1}{c^{2}}\frac{\partial u}{\partial t} - \nabla^{2} u = S(x,t)$& $\left(-\frac{1}{\alpha^{2}}\nabla^{2}+1 \right)\left(u^{n+1}+2u^{n}+u^{n-1}\right) = 4u^{n} + \frac{4}{\alpha^{2}}S(x,t^{n})$  \\
		To  & Diffusive Scheme, $\alpha = \frac{\sqrt{2}}{c \Delta t}$: \\  
		Modified Helmholtz Equation & $\left(-\frac{1}{\alpha^{2}}\nabla^{2}+1 \right)u^{n+1} = \frac{1}{2}\left(5u^{n}-4u^{n-1}+u^{n-2}\right) + \frac{1}{\alpha^{2}} S(x,t^{n+1})$ \\
		
		\hline
		Dimensionally Split  & $\left(-\frac{1}{\alpha^{2}}\nabla^{2}+1 \right)u=f$ $\Rightarrow$ \\
		Modified Helmholtz Equation& $\left(-\frac{1}{\alpha^{2}}\frac{\partial^{2}}{\partial x^{2}}+1 \right)\left(-\frac{1}{\alpha^{2}}\frac{\partial^{2}}{\partial y^{2}}+1 \right)u = f$ $\Rightarrow$ \\
		(2D) &$\left(-\frac{1}{\alpha^{2}}\frac{\partial^{2}}{\partial x^{2}}+1 \right)w = f$, \quad $ \left(-\frac{1}{\alpha^{2}}\frac{\partial^{2}}{\partial y^{2}}+1 \right)u = w$  \\
		\hline
		1D Integral Solution & $\left(-\frac{1}{\alpha^{2}}\frac{d^{2}}{dx^{2}}+1 \right)u = f$ on $(a,b)$ $\Rightarrow$  \\
		& $u(x) = \frac{\alpha}{2}\int_{a}^{b} f(x')e^{-\alpha|x-x'|} \, dx' + Ae^{-\alpha(x-a)} + Be^{-\alpha(b-x)}$ \\
		& $= I[f](x) + Ae^{-\alpha(x-a)} + Be^{-\alpha(b-x)}$\\
		\hline
		1D BC Correction Coefficients & \\
		Dirichlet: & $A = \frac{(u_{a}-I_{a})-\mu (u_{b}-I_{b})}{1-\mu^{2}}$, $B = \frac{(u_{b}-I_{b})-\mu (u_{a}-I_{a})}{1-\mu^{2}}$   \\
		$u(a) = u_{a}$, $u(b) = u_{b}$ & \\
		Neumann: & $A = \frac{\mu(v_{b}+\alpha I_{b})-(v_{a}-\alpha I_{a})}{\alpha \left(1-\mu^2\right)}$, $B = \frac{(v_{b}+\alpha I_{b})-\mu(v_{a}-\alpha I_{a})}{\alpha \left(1-\mu^2\right)}$\\
		$u'(a) = v_{a}$, $u'(b) = v_{b}$ & \\
		Periodic: & $A = \frac{I_{b}}{1-\mu}$, $B = \frac{I_{a}}{1-\mu}$ \\
		$u(a) = u(b)$, $u'(a) = u'(b)$ & \\
		& $I_{a} = I[f](a)$, $I_{b} = I[f](b)$, $\mu = e^{-\alpha(b-a)}$ \\
		\hline
		Fast Convolution Algorithm & $a = x_{0} < x_{1} < \cdots < x_{N} = b$, \\
		&  $x_{j+1} - x_{j}+ \Delta x$, $j=0,...,N-1$ \\
		& $I_{j} = I[f](x_{j}) = \frac{\alpha}{2}\int_{a}^{b} f(x')e^{-\alpha|x_{j}-x'|} \, dx' = I^{L}_{j} + I^{R}_{j}$ \\
		& $I^{L}_{j} = \frac{\alpha}{2}\int_{a}^{x_{j}} f(x')e^{-\alpha|x_{j}-x'|} \, dx'$, $I^{R}_{j} = \frac{\alpha}{2}\int_{x_{j}}^{b} f(x')e^{-\alpha|x_{j}-x'|} \, dx'$ \\
		& $I^{L}_{0}=0$, $I^{L}_{j} = e^{-\alpha \Delta x} I^{L}_{j-1} + J^{L}_{j}$, \\
		& $J^{L}_{j} = \frac{\alpha}{2}\int_{x_{j-1}}^{x_{j}} f(x')e^{-\alpha|x_{j}-x'|} \, dx'$, $j=1,...,N$ \\
		& $I^{R}_{N}=0$, $I^{R}_{j} = e^{-\alpha \Delta x} I^{R}_{j+1} + J^{R}_{j}$, \\
		& $J^{R}_{j} = \frac{\alpha}{2}\int_{x_{j}}^{x_{j+1}} f(x')e^{-\alpha|x_{j}-x'|} \, dx'$, $j=N-1,...,0$ \\
		\hline
		Second Order Quadrature & $J^{R}_{j} = Pf(x_{j})+Qf(x_{j+1}) + R(f(x_{j+1})-2f(x_{j})+f(x_{j-1}))$ \\
		& $J^{L}_{j} = Pf(x_{j})+Qf(x_{j-1}) + R(f(x_{j+1})-2f(x_{j})+f(x_{j-1}))$ \\
		& $P = 1-\frac{1-d}{\nu}$, $Q = -d+\frac{1-d}{\nu}$, $R = \frac{1-d}{\nu^{2}}-\frac{1+d}{2\nu}$ \\
		& $\nu = \alpha \dx$, $d = e^{-\nu}$ \\
		\hline
	\end{tabular}
\end{center}


\section{Treatment of point sources, and soft sources}
\label{sec:sources}

We now consider the inclusion of source terms. We present the algorithm in 1D for simplicity, and observe that the extensions to multiple dimensions are analogous to those shown for dimensional splitting presented in section \ref{sec:molt}.

We are predominantly interested in the case where $S(x,t)$ consists of a large number of time dependent point sources. However, it is often the case that in electromagnetics problems, a soft source is prescribed to excite waves of a  prescribed frequency, or range of frequencies, within the domain. A soft source is so named because, although incident fields are generated at a prescribed fixed spatial location, no scattered fields are generated.

The implementation of a soft source $\sigma(t)$ at $x = x_s$ is accomplished by prescribing the source condition
\begin{equation}
\label{eqn:condition}
u(x_s,t) = \sigma(t).
\end{equation}
However, it can be shown that if we set
\begin{align}
	S(x,t)= \frac{2}{c}\sigma'(t) \delta(x-x_s)
\end{align}
and insert it into the wave equation, then the soft source condition \eqref{eqn:condition} is satisfied, and the solutions are equivalent. Thus, a soft source is nothing more than a point source, whose time-varying field is integrated by the wave equation.

Upon convolving this source term with the Green's function according to \eqref{eqn:Iu}, we find
\begin{align*}
	I\left[\frac{1}{\alpha^2}S\right](x)=& \frac{1}{2\alpha}\int_a^b \left(\frac{2}{c}\sigma'(t_n)\delta(x-x_s) \right)e^{-\alpha|x-y|}dy \\
	=& \frac{\Delta t}{\beta} \sigma'(t_n) e^{-\alpha|x-x_s|},
\end{align*}
where the definition of $\alpha = \beta/(c\Delta t)$ has been utilized.

\begin{remark}
	It is often the case that taking the analytical derivative $\sigma'(t_n)$ is to be avoided, for various reasons. In this case, any finite difference approximation which is of the desired order of accuracy can be substituted.
\end{remark}

Likewise for general point sources,
\[
S(x,t) = \sum_{i} \tilde{\sigma}_i(t) \delta(x-x_i)
\]
the corresponding form of the source term is
\begin{equation}
\label{eqn:Delta_Source}
I\left[\frac{1}{\alpha^2}S\right](x) = \frac{c \Delta t}{2\beta}\sum_i  \tilde{\sigma}_i(t_n) e^{-\alpha|x-x_i|}
\end{equation}
Therefore, it suffices to consider delta functions both for the implementation of soft sources, as well as including time dependent point sources.

\bibliographystyle{amsplain}
\bibliography{ref}

\providecommand{\bysame}{\leavevmode\hbox to3em{\hrulefill}\thinspace}
\providecommand{\MR}{\relax\ifhmode\unskip\space\fi MR }
\providecommand{\MRhref}[2]{%
  \href{http://www.ams.org/mathscinet-getitem?mr=#1}{#2}
}
\providecommand{\href}[2]{#2}
\begin{thebibliography}{10}

\bibitem{Alpert2000}
B.~Alpert, L.~Greengard, and T.~Hagstrom, \emph{{An Integral Evolution Formula
  for the Wave Equation}}, Journal of Computational Physics \textbf{162}
  (2000), no.~2, 536--543.

\bibitem{Alpert2002}
\bysame, \emph{{Nonreflecting Boundary Conditions for the Time-Dependent Wave
  Equation}}, Journal of Computational Physics \textbf{180} (2002), no.~1,
  270--296.

\bibitem{Barnes1986}
J.~{Barnes} and P.~{Hut}, \emph{{A hierarchical O(N log N) force-calculation
  algorithm}}, Nature \textbf{324} (1986), 446--449.

\bibitem{Birdsall1976}
C.K. Birdsall and A.B. Langdon, \emph{Plasma physics via computer simulation},
  Course notes for electrical engineering and computer sciences, McGraw-Hill,
  1976.

\bibitem{birkhoff1906general}
G.~D. Birkhoff, \emph{General mean value and remainder theorems with
  applications to mechanical differentiation and quadrature}, Transactions of
  the American Mathematical Society \textbf{7} (1906), no.~1, pp. 107--136
  (English).

\bibitem{Bruno2010}
O.~P. Bruno and M.~Lyon, \emph{{High-order unconditionally stable FC-AD solvers
  for general smooth domains I . Basic elements}}, Journal of Computational
  Physics \textbf{229} (2010), no.~6, 2009--2033.

\bibitem{Causley2013a}
M.~Causley and A.~Christlieb, \emph{{High order A-stable schemes for the wave
  equation using successive convolution}}, SIAM Journal of Numerical Analysis
  \textbf{52} (2014), no.~1, 220--235.

\bibitem{Causley2013}
M.~Causley, A.~Christlieb, B.~Ong, and L.~Van~Groningen, \emph{{Method of Lines
  Transpose: An Implicit Solution to the Wave Equation}}, Mathematics of
  Computation \textbf{83} (2014), no.~290, 2763--2786.

\bibitem{Causley_Christlieb_Cho}
M.~F. {Causley}, H.~{Cho}, A.~J. {Christlieb}, and D.~C. {Seal}, \emph{{Method
  of lines transpose: High order L-stable O(N) schemes for parabolic equations
  using successive convolution}}, ArXiv e-prints (2015).

\bibitem{Cheng1999}
H.~Cheng, L.~Greengard, and V.~Rokhlin, \emph{{A fast adaptive multipole
  algorithm in three dimensions}}, J. Comput. Phys. \textbf{155} (1999), no.~2,
  468--498.

\bibitem{Cheng2006}
H.~Cheng, J.~Huang, and T.~J. Leiterman, \emph{{An adaptive fast solver for the
  modified Helmholtz equation in two dimensions}}, Journal of Computational
  Physics \textbf{211} (2006), no.~2, 616--637.

\bibitem{Christlieb2004}
A.~J. Christlieb, R.~Krasny, and J.~P. Verboncoeur, \emph{{A treecode algorithm
  for simulating electron dynamics in a Penning-Malmberg trap}}, Computer
  Physics Communications \textbf{164} (2004), no.~1-3, 306--310.

\bibitem{Coifman1993}
R.~Coifman, V.~Rokhlin, and S.~Wandzura, \emph{{The fast multipole method for
  the wave equation: A pedestrian prescription}}, IEEE Trans. Antennas and
  Propagation \textbf{35} (1993), no.~3, 7--12.

\bibitem{Fornberg}
B.~Fornberg, \emph{{A Short Proof of the Unconditional Stability of the
  ADI-FDTD Scheme}},  \textbf{9810751}, 5--8.

\bibitem{Fornberga}
J.~Fornberg, B.and~Zuev and J.~Lee, \emph{{Stability and accuracy of
  time-extrapolated ADI-FDTD methods for solving wave equations}},
  \textbf{9810751}, no.~November 2005.

\bibitem{Gimbutas2002}
Z.~Gimbutas and V.~Rokhlin, \emph{{A generalized fast multipole method for
  nonoscillatory kernels}}, SIAM J. Sci. Comput. \textbf{24} (2002), no.~3,
  796--817.

\bibitem{greengard1987fast}
L.~Greengard and V.~Rokhlin, \emph{{A fast algorithm for particle
  simulations}}, J. Comput. Phys. \textbf{73} (1987), no.~2, 325--348.

\bibitem{kreiss2004difference}
H.~O. Kreiss, N.~A. Petersson, and J.~Ystr{\"o}m, \emph{Difference
  approximations of the neumann problem for the second order wave equation},
  SIAM Journal on Numerical Analysis \textbf{42} (2004), no.~3, 1292--1323.

\bibitem{Li2006}
J.~R. Li, \emph{{Low order approximation of the spherical nonreflecting
  boundary kernel for the wave equation}}, Linear Algebra and its Applications
  \textbf{415} (2006), no.~2-3, 455--468.

\bibitem{Li2009}
P.~Li, H.~Johnston, and R.~Krasny, \emph{{A Cartesian treecode for screened
  coulomb interactions}}, Journal of Computational Physics \textbf{228} (2009),
  no.~10, 3858--3868.

\bibitem{Lindsay2001}
K.~Lindsay and R.~Krasny, \emph{{A Particle Method and Adaptive Treecode for
  Vortex Sheet Motion in Three-Dimensional Flow}}, Journal of Computational
  Physics \textbf{172} (2001), no.~2, 879--907.

\bibitem{Lyon2010}
M.~Lyon and O.~P. Bruno, \emph{{High-order unconditionally stable FC-AD solvers
  for general smooth domains II . Elliptic , parabolic and hyperbolic PDEs ;
  theoretical considerations}}, Journal of Computational Physics \textbf{229}
  (2010), no.~9, 3358--3381.

\bibitem{Namiki2000}
T.~Namiki, \emph{3-d adi-fdtd method-unconditionally stable time-domain
  algorithm for solving full vector maxwell's equations}, Microwave Theory and
  Techniques, IEEE Transactions on \textbf{48} (2000), no.~10, 1743--1748.

\bibitem{Smithe2009}
D.~N. Smithe, J.~R. Cary, and J.~A. Carlsson, \emph{Divergence preservation in
  the adi algorithms for electromagnetics}, J. Comput. Phys. \textbf{228}
  (2009), no.~19, 7289--7299.

\end{thebibliography}

\end{document}